# NEW RESULTS IN TRAJECTORY-BASED SMALL-GAIN WITH APPLICATION TO THE STABILIZATION OF A CHEMOSTAT[§]

Iasson Karafyllis[*] and Zhong-Ping Jiang[**]

[*]Dept. of Environmental Eng., Technical University of Crete,
73100, Chania, Greece, email: ikarafyl@enveng.tuc.gr

[**]Dept. of Electrical and Computer Eng., Polytechnic Institute of New York University,
Six Metrotech Center, Brooklyn, NY 11201, U.S.A., email: zjiang@control.poly.edu

**Abstract**

New trajectory-based small-gain results are obtained for nonlinear feedback systems under relaxed assumptions. Specifically, during a transient period, the solutions of the feedback system may not satisfy some key inequalities that previous small-gain results usually utilize to prove stability properties. The results allow the application of the small-gain perspective to various systems which satisfy less demanding stability notions than the Input-to-Output Stability property. The robust global feedback stabilization problem of an uncertain time-delayed chemostat model is solved by means of the trajectory-based small-gain results.

## 1. Introduction

Small-gain results are important tools for robustness analysis and robust controller design in Mathematical Control Theory. A nonlinear, generalized small-gain theorem was developed in [14], based on the notion of Input-to-State Stability (ISS) originally introduced by Sontag [38]. Recently, nonlinear small-gain results were developed for monotone systems, an important class of nonlinear systems in mathematical biology (see [1,7]). Further extensions of the small-gain perspective to the cases of non-uniform in time stability, discrete-time systems and Lyapunov characterizations are pursued by several authors independently; see, for instance, [2,4,5,9,13,15,16,17,18,19,20,23,41,42]. A general vector small-gain result, which can be applied to a wide class of control systems, was developed in [28].

One of the most important obstacles in applying nonlinear small-gain results is the representation of the original composite system as the feedback interconnection of subsystems which satisfy the Input-to-Output Stability (IOS) property. More specifically, sometimes the subsystems do not satisfy the IOS property: there is a transient period after which the solution enters a certain region of the state space. Within this region of the state space the subsystems satisfy the small-gain requirements. In other words, the essential inequalities, which small-gain results utilize in order to prove stability properties, do not hold for all times: this feature excludes all available small-gain results from possible application. Particularly, this feature is important in systems of Mathematical Biology and Population Dynamics. Indeed, the idea of developing stability results which utilize certain Lyapunov-like conditions after an initial transient period was used in [26,27] with primary motivation from addressing robust feedback stabilization problems for certain chemostat models.

In this work we present small-gain results which can allow a transient period during which the solutions do not satisfy the IOS inequalities (Theorem 2.5 and Theorem 2.6). The obtained results are direct extensions of the recent vector small-gain result in [28] and if the initial transient vanishes then the results coincide with Theorem 3.1 in [28]. The significance of the obtained results is twofold:

- it allows the application of the small-gain perspective to various systems which satisfy less demanding stability notions than ISS,

- it allows the study of systems in Mathematical Biology and Population Dynamics.

---

[§] This work has been supported in part by the NSF grants DMS-0504462 and DMS-0906659.

To emphasize the latter point, we show how the obtained small-gain results can be used for the feedback stabilization of uncertain chemostat models. Chemostat models are often adequately represented by a simple dynamic model involving two state variables, the microbial biomass concentration $X$ and the limiting nutrient concentration $S$ (see [37]). The common delay-free model for microbial growth on a limiting substrate in a chemostat is of the form:

$$\dot{X}(t) = (\mu(S(t)) - D(t))X(t)$$
$$\dot{S}(t) = D(t)(S_i - S(t)) - K\mu(S(t))X(t) \quad (1.1)$$
$$X(t) \in (0, +\infty), S(t) \in (0, S_i), D(t) \geq 0$$

where $S_i$ is the feed substrate concentration, $D$ is the dilution rate (which is used as the control input), $\mu(S)$ is the specific growth rate and $K > 0$ is a biomass yield factor. The literature on control studies of chemostat models of the form (1.1) is extensive. In [6], feedback control of the chemostat by manipulating the dilution rate was studied for the promotion of coexistence. Other interesting control studies of the chemostat can be found in [3,8,10,24,29,30,31]. The stability and robustness of periodic solutions of the chemostat was studied in [32,33]. The problem of the stabilization of a non-trivial steady state $(X_s, S_s)$ of the chemostat model (1.1) was considered in [29], where it was shown that the simple feedback law $D = \mu(S)X / X_s$ is a globally stabilizing feedback. See also the recent work [24] for the study of the robustness properties of the closed-loop system (1.1) with $D = \mu(S)X / X_s$ for time-varying inlet substrate concentration $S_i$. The recent work [26] studied the sampled-data stabilization of the non-trivial steady state $(X_s, S_s)$ of the chemostat model (1.1), while uncertain chemostat models were considered in [27].

In this work we consider the robust global feedback stabilization problem for the chemostat model with delays:

$$\dot{X}(t) = (p(T_r(t)S) - D(t) - b)X(t)$$
$$\dot{S}(t) = D(t)(S_i - S(t)) - K(S(t))\mu(S(t))X(t) \quad (1.2)$$
$$X(t) \in (0, +\infty), S(t) \in (0, S_i), D(t) \geq 0$$

where $T_r(t)S = \{S(t+\theta) : \theta \in [-r, 0]\}$ is the $r$-history of $S$, $b \geq 0$ is the cell mortality rate, $r \geq 0$ is the maximum delay, $K(S) > 0$ is a possibly variable yield coefficient and $p : C^0([-r, 0]; (0, S_i)) \to (0, +\infty)$ is a continuous functional that satisfies

$$\min_{t-r \leq \tau \leq t} \mu(S(\tau)) \leq p(T_r(t)S) \leq \max_{t-r \leq \tau \leq t} \mu(S(\tau)) \quad (1.3)$$

The functions $\mu : [0, S_i] \to [0, \mu_{\max}]$, $K : [0, S_i] \to (0, +\infty)$ with $\mu(0) = 0$, $\mu(S) > 0$ for all $S > 0$, are assumed to be locally Lipschitz functions. The chemostat model (1.2) under (1.3) is very general, since we may have:

- $p(T_r(t)S) = \mu(S(t))$, which gives the standard chemostat model with no delays,

- $p(T_r(t)S) = \mu(S(t-r))$, which gives the time-delayed chemostat model studied in [37],

- $p(T_r(t)S) = \lambda \sum_{i=0}^{n} w_i \mu(S(t-r_i)) + (1-\lambda) \int_{t-r}^{t} h(\tau + r - t)\mu(S(\tau))d\tau$, where $\lambda \in [0,1]$, $h \in C^0([0, r]; [0, +\infty))$

  with $\int_0^r h(s)ds = 1$, $w_i \geq 0$, $r_i \in [0, r]$ ($i = 0, ..., n$) with $\sum_{i=0}^{n} w_i = 1$.

Moreover, it should be noted that the case of variable yield coefficients has been studied recently (see [46,47]) and has been proposed for the justification of experimental results. The reader should notice that chemostat models with time delays were considered in [44,45]. We assume the existence of a non-trivial equilibrium point for (1.2), i.e., the existence of $(S_s, X_s) \in (0, S_i) \times (0, +\infty)$ such that

$$\mu(S_s) = D_s + b \quad , \quad X_s = \frac{D_s(S_i - S_s)}{K(S_s)(D_s + b)} \quad (1.4)$$

where $D_s > 0$ is the equilibrium value for the dilution rate. The stabilization problem for the equilibrium point $(S_s, X_s) \in (0, S_i) \times (0, +\infty)$ is crucial: in [37] it is shown that the equilibrium point is unstable even if $\mu : (0, S_i) \to (0, \mu_{\max}]$ is monotone (e.g., the Monod specific growth rate). Moreover, as remarked in [37] the chemostat model (1.2) under (1.3) allows the expression of the effect of the time difference between consumption of nutrient and growth of the cells (see the discussion on pages 238-240 in [37]). We solve the feedback stabilization problem for the chemostat by providing a *delay-free* feedback which achieves global stabilization (see Theorem 4.1 below). The proof of the theorem relies on the trajectory-based small-gain results of the paper. No knowledge of the maximum delay $r \geq 0$ is assumed.

The structure of the present work is as follows: Section 2 contains the statements of the trajectory-based small-gain results (Theorem 2.5 and Theorem 2.6). Section 3 provides illustrative examples which show the applicability of the obtained results to systems which satisfy less demanding stability notions than ISS. Section 4 is devoted to the development of the solution of the feedback stabilization problem for the uncertain chemostat (1.2). The conclusions are provided in Section 5. The proofs of the trajectory-based small-gain results are given in Appendix A. Finally, for readers' convenience, the definitions of the system-theoretic notions used in this work are given in Appendix B.

**Notations** Throughout this paper, we adopt the following notations:
* We denote by $K^+$ the class of positive, continuous functions defined on $\Re_+ := \{x \in \Re : x \geq 0\}$. We say that a function $\rho : \Re_+ \to \Re_+$ is positive definite if $\rho(0) = 0$ and $\rho(s) > 0$ for all $s > 0$. By $K$ we denote the set of positive definite, increasing and continuous functions. We say that a positive definite, increasing and continuous function $\rho : \Re_+ \to \Re_+$ is of class $K_\infty$ if $\lim_{s \to +\infty} \rho(s) = +\infty$. By $KL$ we denote the set of all continuous functions $\sigma = \sigma(s, t) : \Re_+ \times \Re_+ \to \Re_+$ with the properties: (i) for each $t \geq 0$ the mapping $\sigma(\cdot, t)$ is of class $K$; (ii) for each $s \geq 0$, the mapping $\sigma(s, \cdot)$ is non-increasing with $\lim_{t \to +\infty} \sigma(s, t) = 0$.

* By $\| \ \|_X$, we denote the norm of the normed linear space $X$. By $| \ |$ we denote the Euclidean norm of $\Re^n$. Let $U \subseteq X$ with $0 \in U$. By $B_U[0, r] := \{ u \in U ; \|u\|_X \leq r \}$ we denote the intersection of $U \subseteq X$ with the closed ball of radius $r \geq 0$, centered at $0 \in U$. If $U \subseteq \Re^n$ then $\text{int}(U)$ denotes the interior of the set $U \subseteq \Re^n$.

* $x'$ denotes the transpose of $x$.

* $\Re_+^n := (\Re_+)^n = \{(x_1, ..., x_n)' \in \Re^n : x_1 \geq 0, ..., x_n \geq 0\}$. $\{e_i\}_{i=1}^n$ denotes the standard basis of $\Re^n$. $Z_+$ denotes the set of non-negative integers.

* Let $x, y \in \Re^n$. We say that $x \leq y$ if and only if $(y - x) \in \Re_+^n$. We say that a function $\rho : \Re_+^n \to \Re^+$ is of class $N_n$, if $\rho$ is continuous with $\rho(0) = 0$ and such that $\rho(x) \leq \rho(y)$ for all $x, y \in \Re_+^n$ with $x \leq y$.

* For $t \geq t_0 \geq 0$ let $[t_0, t] \ni \tau \to V(\tau) = (V_1(\tau), ..., V_n(\tau))' \in \Re^n$ be a bounded map. We define $[V]_{[t_0, t]} := \left( \sup_{\tau \in [t_0, t]} V_1(\tau), ..., \sup_{\tau \in [t_0, t]} V_n(\tau) \right)$. For a measurable and essentially bounded function $x : [a, b] \to \Re^n$ $ess \sup_{t \in [a,b]} |x(t)|$ denotes the essential supremum of $|x(\cdot)|$. Given a function $x : [a - r, b) \to \Re^n$, where $r > 0$, $a < b$, we define $T_r(t)x := x(t + \theta) ; \theta \in [-r, 0]$, for $t \in [a, b)$, to be the $r$-history of $x$.

* We say that $\Gamma : \Re_+^n \to \Re_+^m$ is non-decreasing if $\Gamma(x) \leq \Gamma(y)$ for all $x, y \in \Re_+^n$ with $x \leq y$. For an integer $k \geq 1$, we define $\Gamma^{(k)}(x) = \underbrace{\Gamma \circ \Gamma \circ ... \circ \Gamma}_{k \ times}(x)$, when $m = n$.

* We define $\mathbf{1} = (1, 1, ..., 1)' \in \Re^n$. If $u, v \in \Re$ and $u \leq v$ then $\mathbf{1}u \leq \mathbf{1}v$.

* Let $U$ be a subset of a normed linear space $\mathcal{U}$, with $0 \in U$. By $M(U)$ we denote the set of all locally bounded functions $u : \Re_+ \to U$. By $u_0$ we denote the identically zero input, i.e., the input that satisfies $u_0(t) = 0 \in U$ for all $t \geq 0$. If $U \subseteq \Re^n$ then $M_U$ denotes the space of measurable, locally bounded functions $u : \Re_+ \to U$.

* Let $A \subseteq X$, $B \subseteq Y$, where $X, Y$ are normed linear spaces. We denote by $C^0(A; B)$ the class of continuous mappings $f : A \to B$. For $x \in C^0([-r, 0]; \Re^n)$ we define $\|x\|_r := \max_{\theta \in [-r, 0]} |x(\theta)|$.

## 2. New Trajectory-Based Small-Gain Theorems

In this Section we state the main results of the present work. The proofs of the main results (Theorem 2.5 and Theorem 2.6) are provided at Appendix A. For the statement of the main result one needs to know the abstract system theoretic framework introduced in [21,22,23] and used in [28]. For the convenience of the reader, all definitions of the basic notions are provided in Appendix B.

The following technical definitions were used in [28] and are needed here.

**Definition 2.1:** Let $x = (x_1,...,x_n)' \in \Re^n$, $y = (y_1,...,y_n)' \in \Re^n$. We define $z = MAX\{x, y\}$, where $z = (z_1,...,z_n)' \in \Re^n$ satisfies $z_i = \max\{x_i, y_i\}$ for $i = 1,...,n$. Similarly for $u_1,...,u_m \in \Re^n$ we have $z = MAX\{u_1,...,u_m\}$ is a vector $z = (z_1,...,z_n)' \in \Re^n$ with $z_i = \max\{u_{1i},...,u_{mi}\}$, $i = 1,...,n$.

**Definition 2.2:** We say that $\Gamma : \Re_+^n \to \Re_+^n$ is **MAX-preserving** if $\Gamma : \Re_+^n \to \Re_+^n$ is non-decreasing and for every $x, y \in \Re_+^n$ the following equality holds:
$$\Gamma(MAX\{x, y\}) = MAX\{\Gamma(x), \Gamma(y)\} \tag{2.1}$$

The above defined MAX-preserving maps enjoy the following important property (see [28]).

**Proposition 2.3:** $\Gamma : \Re_+^n \to \Re_+^n$ with $\Gamma(x) = (\Gamma_1(x),...,\Gamma_n(x))'$ is MAX-preserving if and only if there exist non-decreasing functions $\gamma_{i,j} : \Re_+ \to \Re_+$, $i, j = 1,...,n$ with $\Gamma_i(x) = \max_{j=1,...,n} \gamma_{i,j}(x_j)$ for all $x \in \Re_+^n$, $i = 1,...,n$.

The following class of MAX-preserving mappings plays an important role in what follows.

**Definition 2.4:** Let $\Gamma : \Re_+^n \to \Re_+^n$ with $\Gamma(x) = (\Gamma_1(x),...,\Gamma_n(x))'$ be a MAX-preserving mapping for which there exist functions $\gamma_{i,j} \in N_1$, $i, j = 1,...,n$ with $\Gamma_i(x) = \max_{j=1,...,n} \gamma_{i,j}(x_j)$ for all $x \in \Re_+^n$, $i = 1,...,n$. We say that $\Gamma : \Re_+^n \to \Re_+^n$ satisfies the cyclic small-gain conditions if the following inequalities hold:
$$\gamma_{i,i}(s) < s, \ \forall s > 0, \ i = 1,...,n \tag{2.2}$$
and if $n > 1$ then for each $r = 2,...,n$ it holds that:
$$\left(\gamma_{i_1,i_2} \circ \gamma_{i_2,i_3} \circ ... \circ \gamma_{i_r,i_1}\right)(s) < s, \ \forall s > 0 \tag{2.3}$$
for all $i_j \in \{1,...,n\}$, $i_j \neq i_k$ if $j \neq k$.

Proposition 2.7 in [28] shows that the MAX-preserving continuous mapping $\Gamma : \Re_+^n \to \Re_+^n$ satisfies the cyclic small-gain conditions if and only if $0 \in \Re^n$ is Globally Asymptotically Stable for the discrete-time $x(k+1) = \Gamma(x(k))$, where $x(k) \in \Re_+^n$, $k \in Z^+$. The following facts are consequences of the related results in [5,28,36,43] and definitions (2.1), (2.2), (2.4) and will be used repeatedly in the proofs of the main results of the present section.

**Fact I:** If $\Gamma : \Re_+^n \to \Re_+^n$ satisfies the cyclic small-gain conditions, then $\lim_{k \to +\infty} \Gamma^{(k)}(x) = 0$ for all $x \in \Re_+^n$ and $\Gamma^{(k)}(x) \leq Q(x) = MAX\{x, \Gamma(x), \Gamma^{(2)}(x),..., \Gamma^{(n-1)}(x)\}$ for all $k \geq 1$ and $x \in \Re_+^n$.

**Fact II:** If $\Gamma : \Re_+^n \to \Re_+^n$ is a MAX-preserving mapping, then the mapping $Q(x) = MAX\{x, \Gamma(x), \Gamma^{(2)}(x),..., \Gamma^{(n-1)}(x)\}$ is a MAX-preserving mapping.

**Fact III:** If $\Gamma : \Re_+^n \to \Re_+^n$ satisfies the cyclic small-gain conditions, then $\Gamma(Q(x)) \leq Q(x)$ and $Q(x) \geq x$ for all $x \in \Re_+^n$, where $Q(x) = MAX\{x, \Gamma(x), \Gamma^{(2)}(x),..., \Gamma^{(n-1)}(x)\}$.

**Fact IV:** If $p \in N_n$ and $R: \Re_+^n \to \Re_+^n$ is a non-decreasing mapping, then the following inequality holds for all $s, r \in \Re_+$: $p(MAX\{R(\mathbf{1}s), R(\mathbf{1}r)\}) = \max(p(R(\mathbf{1}s)), p(R(\mathbf{1}r)))$.

**Fact V:** If $\Gamma: \Re_+^n \to \Re_+^n$ satisfies the cyclic small-gain conditions and $x, y \in \Re_+^n$ satisfy $x \leq MAX\{y, \Gamma(x)\}$, then $x \leq Q(y)$, where $Q(x) = MAX\{x, \Gamma(x), \Gamma^{(2)}(x), ..., \Gamma^{(n-1)}(x)\}$.

We consider an abstract control system $\Sigma := (X, Y, M_U, M_D, \phi, \pi, H)$ with the BIC property for which $0 \in X$ is a robust equilibrium point from the input $u \in M_U$ (see Appendix B for the notions of an abstract control system, the BIC property and the notion of a robust equilibrium point). We suppose that there exists a set-valued map $\Re_+ \ni t \to S(t) \subseteq X$ with $0 \in S(t)$ for all $t \geq 0$, maps $V_i : \bigcup_{t \geq 0} \{t\} \times S(t) \times U \to \Re^+$, with $V_i(t,0,0) = 0$ for all $t \geq 0$ ($i = 1,...,n$) and a MAX-preserving continuous map $\Gamma: \Re_+^n \to \Re_+^n$ with $\Gamma(0) = 0$ such that the following hypotheses hold:

**(H1)** There exist functions $\sigma \in KL$, $\zeta \in N_1$, $L: \bigcup_{t \geq 0} \{t\} \times S(t) \to \Re^+$ with $L(t,0) = 0$ for all $t \geq 0$, such that for every $(t_0, x_0, u, d) \in \Re_+ \times X \times M_U \times M_D$ with $\phi(t, t_0, x_0, u, d) \in S(t)$ for all $t \in [t_0, t_{max})$ the mappings $t \to V(t) = (V_1(t, \phi(t, t_0, x_0, u, d), u(t)), ..., V_n(t, \phi(t, t_0, x_0, u, d), u(t)))'$ and $t \to L(t) = L(t, \phi(t, t_0, x_0, u, d))$ are locally bounded on $[t_0, t_{max})$ and the following estimates hold:

$$V(t) \leq MAX\left\{\mathbf{1}\sigma(L(t_0), t - t_0), \Gamma([V]_{[t_0,t]}), \mathbf{1}\zeta\left(\|u\|_{U_{[t_0,t]}}\right)\right\}, \text{ for all } t \in [t_0, t_{max}) \quad (2.4)$$

where $t_{max}$ is the maximal existence time of the transition map of $\Sigma$.

**(H2)** For every $(t_0, x_0, u, d) \in \Re_+ \times X \times M_U \times M_D$ there exists $\xi \in \pi(t_0, x_0, u, d)$ such that $\phi(t, t_0, x_0, u, d) \in S(t)$ for all $t \in [\xi, t_{max})$. Moreover, there exist functions $v, c, \tilde{c} \in K^+$, $a, \eta, \tilde{\eta}, p^u, g^u \in N_1$, $p \in N_n$, such that the following inequalities hold for every $(t_0, x_0, u, d) \in \Re_+ \times X \times M_U \times M_D$:

$$L(t) \leq \max\left\{v(t - t_0), c(t_0), a(\|x_0\|_X), p([V]_{[\xi,t]}), p^u\left(\|u\|_{U_{[t_0,t]}}\right)\right\}, \text{ for all } t \in [\xi, t_{max}) \quad (2.5)$$

$$\|\phi(t, t_0, x_0, u, d)\|_X \leq \max\left\{v(t - t_0), \tilde{c}(t_0), a(\|x_0\|_X), \tilde{\eta}\left(\|u\|_{U_{[t_0,t]}}\right)\right\}, \text{ for all } t \in [t_0, \xi] \quad (2.6)$$

$$\xi \leq t_0 + a(\|x_0\|_X) + c(t_0) \quad (2.7)$$

$$\|H(t, \phi(t, t_0, x_0, u, d), u(t))\|_Y \leq \max\left\{a(c(t_0)\|x_0\|_X), \eta\left(\|u\|_{U_{[t_0,t]}}\right)\right\}, \text{ for all } t \in [t_0, \xi] \quad (2.8)$$

$$L(\xi, \phi(\xi, t_0, x_0, u, d)) \leq \max\left\{a(c(t_0)\|x_0\|_X), g^u\left(\|u\|_{U_{[t_0,\xi]}}\right)\right\} \quad (2.9)$$

**(H3)** There exist functions $b \in N_1$, $g \in N_n$, $\mu, \kappa \in K^+$ such that the following inequalities hold:

$$\mu(t)\|x\|_X \leq b(L(t,x) + g(V(t,x,u)) + \kappa(t)), \text{ for all } (t, x, u) \in \bigcup_{t \geq 0} \{t\} \times S(t) \times U \quad (2.10)$$

where $V(t, x, u) = (V_1(t, x, u), ..., V_n(t, x, u))'$.

**(H4)** *There exists $q \in N_n$ such that the following inequality holds*:

$$\|H(t,x,u)\|_Y \leq q\big(V(t,x,u)\big), \text{ for all } (t,x,u) \in \bigcup_{t \geq 0}\{t\} \times S(t) \times U \tag{2.11}$$

**Discussion of Hypotheses (H1), (H2):** Hypotheses (H1), (H2) hold automatically when hypotheses (H1-3) of Theorem 3.1 in [28] hold (hypotheses (H1-3) in [28] correspond to the special case $S(t) = X$ and $\xi = t_0$). Consequently, Hypotheses (H1), (H2) are less restrictive hypotheses. Indeed, inequalities (2.4), (2.5) are not assumed to hold for all times $t \in [t_0, t_{max})$ but only after the solution map $\phi(t, t_0, x_0, u, d)$ has entered the set $S(t) \subseteq X$. Moreover, the set-valued map $S(t) \subseteq X$ is not assumed to be positively invariant.

We are now ready to state the main results.

**Theorem 2.5 (Trajectory-Based Small-Gain Result for IOS):** *Consider system $\Sigma := (X, Y, M_U, M_D, \phi, \pi, H)$ under the above hypotheses. Assume that the MAX-preserving continuous map $\Gamma : \mathfrak{R}_+^n \to \mathfrak{R}_+^n$ with $\Gamma(0) = 0$ satisfies the cyclic small-gain conditions. Then system $\Sigma$ satisfies the IOS property from the input $u \in M_U$ with gain $\gamma(s) := \max\{\eta(s), q(G(s))\}$, where $G(s) = (G_1(s), ..., G_n(s))'$ is defined by:*

$$G(s) = Q\big(\mathbf{1}\max\big\{\sigma\big(p^u(s),0\big), \sigma\big(g^u(s),0\big), \sigma\big(p(Q(\mathbf{1}\sigma(g^u(s),0))),0\big), \sigma\big(p(Q(\mathbf{1}\zeta(s))),0\big), \zeta(s)\big\}\big) \tag{2.12}$$

*with $Q(x) = MAX\{x, \Gamma(x), \Gamma^{(2)}(x), ..., \Gamma^{(n-1)}(x)\}$ for all $x \in \mathfrak{R}_+^n$. Moreover, if $c \in K^+$ is bounded, then system $\Sigma$ satisfies the UIOS property from the input $u \in M_U$ with gain $\gamma(s) := \max\{\eta(s), q(G(s))\}$.*

**Remark:** It is of interest to note that Theorem 2.5 is a new trajectory-based small-gain result for IOS because inequalities (2.4), (2.5) are not assumed to hold for all times. Instead, we assume that for each trajectory there exists a time $\xi \in \pi(t_0, x_0, u, d)$ after which inequalities (2.4), (2.5) hold. On the other hand, in order to be able to conclude IOS for the system, we have to assume additional inequalities which hold for the transient period $t \in [t_0, \xi]$, i.e., inequalities (2.6), (2.7), (2.8), (2.9) are required to hold.

We consider next an abstract control system $\Sigma := (X, Y, M_U, M_D, \phi, \pi, H)$ with $U = \{0\}$ and the BIC property for which $0 \in X$ is a robust equilibrium point from the input $u \in M_U$. Suppose that there exists a set-valued map $\mathfrak{R}_+ \ni t \to S(t) \subseteq X$ with $0 \in S(t)$ for all $t \geq 0$, maps $V_i : \bigcup_{t \geq 0}\{t\} \times S(t) \to \mathfrak{R}^+$, with $V_i(t, 0) = 0$ for all $t \geq 0$ ($i = 1,...,n$) and a MAX-preserving continuous map $\Gamma : \mathfrak{R}_+^n \to \mathfrak{R}_+^n$ with $\Gamma(0) = 0$ such that the following hypothesis holds:

**(H5)** *There exist functions $\sigma \in KL$, $L : \bigcup_{t \geq 0}\{t\} \times S(t) \to \mathfrak{R}^+$ with $L(t,0) = 0$ for all $t \geq 0$, such that for every $(t_0, x_0, d) \in \mathfrak{R}_+ \times X \times M_D$ with $\phi(t, t_0, x_0, u_0, d) \in S(t)$ for all $t \in [t_0, t_{max})$ the mappings $t \to V(t) = (V_1(t, \phi(t, t_0, x_0, u_0, d)), ..., V_n(t, \phi(t, t_0, x_0, u_0, d)))'$ and $t \to L(t) = L(t, \phi(t, t_0, x_0, u_0, d))$ are locally bounded on $[t_0, t_{max})$ and the following estimates hold:*

$$V(t) \leq MAX\big\{\mathbf{1}\sigma\big(L(t_0), t - t_0\big), \Gamma\big([V]_{[t_0,t]}\big)\big\}, \text{ for all } t \in [t_0, t_{max}) \tag{2.13}$$

*where $t_{max}$ is the maximal existence time of the transition map of $\Sigma$.*

**(H6)** *For every $(t_0, x_0, d) \in \mathfrak{R}_+ \times X \times M_D$ there exists $\xi \in \pi(t_0, x_0, u_0, d)$ such that $\phi(t, t_0, x_0, u_0, d) \in S(t)$ for all $t \in [\xi, t_{max})$. Moreover, there exist functions $v, c \in K^+$, $a \in N_1$, $p \in N_n$, such that for every $(t_0, x_0, d) \in \mathfrak{R}_+ \times X \times M_D$ the following inequalities hold:*

$$L(t) \leq \max\{v(t-t_0), c(t_0), a(\|x_0\|_\mathcal{X}), p([V]_{[\xi,t]})\}, \text{ for all } t \in [\xi, t_{\max}) \tag{2.14}$$

$$\|\phi(t,t_0,x_0,u_0,d)\|_\mathcal{X} \leq \max\{v(t-t_0), c(t_0), a(\|x_0\|_\mathcal{X})\}, \text{ for all } t \in [t_0, \xi] \tag{2.15}$$

$$\xi \leq t_0 + a(\|x_0\|_\mathcal{X}) + c(t_0) \tag{2.16}$$

$$L(\xi) \leq a(\|x_0\|_\mathcal{X}) + c(t_0) \tag{2.17}$$

**Discussion of Hypothesis (H6):** Hypothesis (H6) is almost the same with Hypothesis (H2) applied to the case $U = \{0\}$. Nonetheless, notice the difference that the estimate for $L(\xi)$ in inequality (2.17) is less tight than the estimate needed in inequality (2.9) of Hypothesis (H2). Indeed, when $x_0 = 0$, estimate (2.17) does not yield $L(\xi) = 0$, contrary to the estimate (2.9), which gives $L(\xi) = 0$. Finally, the analogue of inequality (2.8) for $U = \{0\}$ is not needed in hypothesis (H6).

**Theorem 2.6 (Trajectory-Based Small-Gain for Robust Global Asymptotic Output Stability (RGAOS)):** Consider system $\Sigma := (\mathcal{X}, \mathcal{Y}, M_U, M_D, \phi, \pi, H)$ with $U = \{0\}$ under hypotheses (H3-6). Assume that the MAX-preserving continuous map $\Gamma : \mathfrak{R}_+^n \to \mathfrak{R}_+^n$ with $\Gamma(0) = 0$ satisfies the cyclic small-gain conditions. Then system $\Sigma$ is RGAOS. Moreover, if $\Sigma := (\mathcal{X}, \mathcal{Y}, M_U, M_D, \phi, \pi, H)$ is $T-$ periodic for certain $T > 0$ then system $\Sigma$ is Uniformly RGAOS (URGAOS).

It is clear that Hypotheses (H5), (H6) is less demanding than Hypotheses (H1), (H2). On the other hand the conclusion of Theorem 2.6 is weaker than the conclusion of Theorem 2.5: Theorem 2.6 guarantees RGAOS while Theorem 2.5 guarantees IOS.

The proofs of Theorem 2.5 and Theorem 2.6 are provided at Appendix A and are similar in spirit to the proof of Theorem 3.1 in [28].

## 3. Examples and Discussions

The first example indicates that the trajectory-based small-gain results of the previous section can be used to study the feedback interconnection of systems which do not necessarily satisfy the IOS property.

**Example 3.1:** Consider the system

$$\begin{aligned} \dot{x} &= f(d, x, y) \\ \dot{y} &= g(d, x, y) \\ x &\in \mathfrak{R}^n, y \in \mathfrak{R}^k, d \in D \subset \mathfrak{R}^l \end{aligned} \tag{3.1}$$

where $D \subset \mathfrak{R}^l$ is a non-empty compact set, $f : D \times \mathfrak{R}^n \times \mathfrak{R}^k \to \mathfrak{R}^n$, $g : D \times \mathfrak{R}^n \times \mathfrak{R}^k \to \mathfrak{R}^k$ are locally Lipschitz mappings with $f(d,0,0) = 0$, $g(d,0,0) = 0$ for all $d \in D$. Suppose that there exist positive definite, continuously differentiable and radially unbounded functions $V : \mathfrak{R}^n \to \mathfrak{R}_+$, $W : \mathfrak{R}^n \to \mathfrak{R}_+$ satisfying the following inequalities for all $(x, y) \in \mathfrak{R}^n \times \mathfrak{R}^k$:

$$\max_{d \in D} \nabla V(x) f(d, x, y) \leq -2 \frac{V(x)}{1+V(x)} + \frac{W(y)}{(1+V(x))(1+W(y))} \tag{3.2}$$

$$\max_{d \in D} \nabla W(y) g(d, x, y) \leq -2 \frac{W(y)}{1+W(y)} + V(x) \tag{3.3}$$

It is clear that the subsystem $\dot{y} = g(d, x, y)$ does not satisfy necessarily the ISS property from the input $x \in \Re^n$. Consequently, the classical small-gain theorem in [14] cannot be applied because the $y-$ subsystem in (3.1) is not ISS but integral ISS with $x \in \Re^n$ as input. Recent small-gain approaches have been used for system (3.1), where it is shown that $0 \in \Re^n \times \Re^k$ is Globally Asymptotically Stable (see [2,12]) for the disturbance-free case. Here we will show, by making use of Theorem 2.6 that $0 \in \Re^n \times \Re^k$ is Uniformly Robustly Globally Asymptotically Stable (URGAS).

Indeed, the inequalities (3.2), (3.3) imply that system (3.1) is RFC. Notice that inequality (3.2) implies for all $(x, y) \in \Re^n \times \Re^k$:

$$\max_{d \in D} \nabla V(x) f(d, x, y) \leq 1$$

The above differential inequality implies that the solution $(x(t), y(t)) \in \Re^n \times \Re^k$ of (3.1) with initial condition $x(0) = x_0 \in \Re^n$, $y(0) = y_0 \in \Re^k$ corresponding to arbitrary $d \in M_D$ satisfies the following inequality for all $t \in [0, t_{\max})$:

$$V(x(t)) \leq V(x_0) + t \quad (3.4)$$

Moreover, inequality (3.3) implies for all $(x, y) \in \Re^n \times \Re^k$:

$$\max_{d \in D} \nabla W(y) g(d, x, y) \leq V(x)$$

The above differential inequality in conjunction with (3.4) implies that the solution $(x(t), y(t)) \in \Re^n \times \Re^k$ of (3.1) with initial condition $x(0) = x_0 \in \Re^n$, $y(0) = y_0 \in \Re^k$ corresponding to arbitrary $d \in M_D$ satisfies the following inequality for all $t \in [0, t_{\max})$:

$$W(y(t)) \leq W(y_0) + tV(x_0) + \frac{t^2}{2} \quad (3.5)$$

Inequalities (3.4), (3.5) imply that system (3.1) is Robustly Forward Complete (see Appendix B). Inequalities (3.2), (3.3) imply for every $\varepsilon > 0$ the existence of a positive definite function $\rho \in C^0(\Re_+; \Re_+)$ such that:

$$V(x) \geq \frac{(1+\varepsilon)W(y)}{2(1+W(y))} \Rightarrow \max_{d \in D} \nabla V(x) f(d, x, y) \leq -\rho(V(x)) \quad (3.6)$$

$$W(y) \geq \frac{(1+\varepsilon)V(x)}{2-V(x)} \text{ and } V(x) < 2 \Rightarrow \max_{d \in D} \nabla W(y) g(d, x, y) \leq -\rho(W(y)) \quad (3.7)$$

Lemma 3.5 in [25] in conjunction with implications (3.6), (3.7) implies that the existence of $\sigma \in KL$ such that the following inequalities hold:

$$V(x(t)) \leq \max\left\{ \sigma(V(x_0), t), \sup_{0 \leq \tau \leq t} \frac{1+\varepsilon}{2} \frac{W(y(\tau))}{1+W(y(\tau))} \right\}, \text{ for all } t \geq 0 \quad (3.8)$$

$$V(x(t)) \leq \max\left\{ \sigma(V(x(\xi)), t-\xi), \sup_{\xi \leq \tau \leq t} \frac{1+\varepsilon}{2} \frac{W(y(\tau))}{1+W(y(\tau))} \right\}, \text{ for all } 0 \leq \xi \leq t \quad (3.9)$$

$$W(y(t)) \leq \max\left\{ \sigma(W(y(\xi)), t-\xi), \sup_{\xi \leq \tau \leq t} (1+\varepsilon) \frac{V(x(\tau))}{2-V(x(\tau))} \right\}, \text{ for all } 0 \leq \xi \leq t \text{ with } \sup_{\xi \leq \tau \leq t} V(x(\tau)) < 2 \quad (3.10)$$

Inequality (3.8) shows that for every $x(0) = x_0 \in \Re^n$, $y(0) = y_0 \in \Re^k$, $d \in M_D$ there exists $\xi \geq 0$ such that $(x(t), y(t)) \in S$ for $t \geq \xi$ with $S := \{(x, y) \in \Re^n \times \Re^k : V(x) \leq 1+\varepsilon\}$. This follows from Proposition 7 in [39] which implies the existence of $a_1, a_2 \in K_\infty$ such that $\sigma(s,t) \leq a_1(\exp(-t)a_2(s))$ for all $s, t \geq 0$. Let $a_3 \in K_\infty$ be a function satisfying $V(x) \leq a_3(|x|)$ for all $x \in \Re^n$. The reader can verify that $(x(t), y(t)) \in S$ for $t \geq \xi$ with

$$\xi := \ln\left(1 + \frac{a_2(a_3(|x_0|))}{a_1^{-1}(1+\varepsilon)}\right) \tag{3.11}$$

Moreover, inequality (2.16) holds for appropriate $a \in K_\infty$ and $c(t) \equiv 1$. Define $V_1(x,y) = \frac{V(x)}{2 - V(x)}$, $V_2(x,y) = W(y)$. The reader should notice that for $\varepsilon \in (0,1)$ inequalities (3.4), (3.5), (3.9), (3.10) and definition (3.11) guarantee that inequalities (2.13), (2.14), (2.15) and (2.17) hold for appropriate $\sigma \in KL$, $\nu \in K^+$, $a \in K_\infty$ with $c(t) \equiv 1$, $p \equiv 0$, $\gamma_{1,2}(s) := \frac{(1+\varepsilon)s}{4 + (3-\varepsilon)s}$, $\gamma_{2,1}(s) := (1+\varepsilon)s$, $\gamma_{1,1}(s) = \gamma_{2,2}(s) \equiv 0$, $L(x,y) := V(x) + W(y)$ and $H(t,x,y) := |(x,y)|$. Finally, notice that the MAX-preserving mapping $\Gamma : \Re_+^2 \to \Re_+^2$ with $\Gamma_i(x) = \max_{j=1,2} \gamma_{i,j}(x_j)$ ($i = 1,2$) satisfies the cyclic small-gain conditions for $\varepsilon \in (0,1)$. The reader can verify that hypotheses (H3), (H4) are automatically satisfied for appropriate functions $b \in N_1$, $g, q \in N_2$, $\mu, \kappa \in K^+$.

By virtue of Theorem 2.6, we conclude that the autonomous system (3.1) is URGAS. ◁

The following example deals with the robust global sampled-data stabilization of a nonlinear planar system.

**Example 3.2:** Consider the following planar system

$$\begin{aligned} \dot{x} &= -(1+y^2)x + y \\ \dot{y} &= f(x) + g(x)y + u \\ (x, y) &\in \Re^2, u \in \Re \end{aligned} \tag{3.12}$$

where $f, g : \Re \to \Re$ are locally Lipschitz functions with $f(0) = 0$. We will show that there exist constants $M > 0$ sufficiently large and $r > 0$ sufficiently small so that system (3.12) in closed loop with the feedback law $u = -My$ applied with zero order hold, i.e., the closed-loop system

$$\begin{aligned} \dot{x}(t) &= -(1+y^2(t))x(t) + y(t) \\ \dot{y}(t) &= -My(\tau_i) + f(x(t)) + g(x(t))y(t)), \, t \in [\tau_i, \tau_{i+1}) \\ \tau_{i+1} &= \tau_i + \exp(-w(\tau_i))r, \quad w(t) \in \Re^+ \end{aligned} \tag{3.13}$$

satisfies the UISS property with zero gain when $w$ is considered as input.

First, notice that there exists a function $\sigma \in KL$ such that for all $(x_0, y) \in \Re^n \times L^\infty_{loc}(\Re^+; \Re)$ the solution of $\dot{x} = -(1 + y^2)x + y$ with initial condition $x(0) = x_0$ corresponding to inputs $y \in L^\infty_{loc}(\Re^+; \Re)$ satisfies the following estimate for all $t \geq 0$:

$$|x(t)| \leq \max\left\{\sigma(|x_0|, t), \sup_{0 \leq \tau \leq t} \gamma(|y(\tau)|)\right\} \tag{3.14}$$

with $\gamma(s) := \frac{\sqrt{2} s}{\sqrt{1+4s^2}}$. Indeed, inequality (3.14) can be verified by using the Lyapunov function $V(x) = x^2$ which satisfies the following implication:

$$\text{if } V(x) = x^2 \geq \frac{2y^2}{1+4y^2} \text{ then } \dot{V} \leq -\frac{1}{4}V(x)$$

The above implication in conjunction with Lemma 3.5 in [25] guarantees that (3.14) holds for appropriate $\sigma \in KL$. Next we show the following claim.

**Claim 1:** *For every $\varepsilon, a > 0$ there exist $\sigma \in KL$, $M > 0$ sufficiently large and $r > 0$ sufficiently small such that for every $(y_0, x, w) \in \Re \times L_{loc}^\infty(\Re^+; B[0,a]) \times L_{loc}^\infty(\Re^+; \Re^+)$ the solution of*

$$\dot{y}(t) = -My(\tau_i) + f(x(t)) + g(x(t))y(t), \; t \in [\tau_i, \tau_{i+1})$$
$$\tau_{i+1} = \tau_i + \exp(-w(\tau_i))\, r \quad, \quad w(t) \in \Re^+ \tag{3.15}$$

*with initial condition $y(0) = y_0$ corresponding to inputs $(x, w) \in L_{loc}^\infty(\Re^+; B[0,a]) \times L_{loc}^\infty(\Re^+; \Re^+)$ satisfies the following inequality:*

$$|y(t)| \leq \max\left\{ \sigma(|y_0|, t),\; \varepsilon \sup_{0 \leq \tau \leq t} |x(\tau)| \right\} \tag{3.16}$$

**Proof of Claim 1:** Let $\varepsilon, a > 0$ be arbitrary. Since $f, g : \Re \to \Re$ are locally Lipschitz functions with $f(0) = 0$, there exist constants $P, Q > 0$ such that

$$|f(x)| \leq P|x| \text{ and } g(x) \leq Q, \text{ for all } x \in B[0,a] \tag{3.17}$$

Let $M > 0$ and $r > 0$ be chosen so that:

$$M \geq 2 + 2Q + \frac{9P^2}{2\varepsilon^2} \text{ and } 3(M+Q)r\exp(Qr) \leq 1 \tag{3.18}$$

Consider a solution $y(t)$ of (3.15) corresponding to arbitrary $(x, w) \in L_{loc}^\infty(\Re^+; B[0,a]) \times L_{loc}^\infty(\Re^+; \Re^+)$ with initial condition $y(0) = y_0 \in \Re$. By virtue of Proposition 2.5 in [21], there exists a maximal existence time for the solution denoted by $t_{\max} \leq +\infty$. Moreover, let $\pi := \{\tau_0, \tau_1, ...\}$ the set of sampling times (which may be finite if $t_{\max} < +\infty$) and $mp(t) := \max\{\tau \in \pi : \tau \leq t\}$. Let $\|x\| := \sup_{0 \leq s \leq t} |x(s)|$ and $\tau = mp(t)$. Inequalities (3.17), (3.18) and the fact that $t - \tau \leq r$ in conjunction with the Gronwall-Bellman inequality implies:

$$|y(t) - y(\tau)| \leq \frac{(M+Q)r\exp(Qr)}{1 - (M+Q)r\exp(Qr)} |y(t)| + \frac{Pr}{1 - (M+Q)r\exp(Qr)} \exp(Qr)\|x\| \tag{3.19}$$

Define $V(t) = y^2(t)$ on $[0, t_{\max})$. Let $I \subset [0, t_{\max})$ be the zero Lebesgue measure set where $y(t)$ is not differentiable or where $\dot{y}(t) \neq -My(\tau_i) + f(x(t)) + g(x(t))y(t))$. Using (3.17), (3.18) and (3.19) we obtain for all $t \in [0, t_{\max}) \setminus I$:

$$\dot{V} \leq -2V(t) + \frac{\varepsilon^2}{2}\|x\|^2, \text{ for all } t \in [0, t_{\max}) \setminus I \tag{3.20}$$

Direct integration of the differential inequality (3.20) and the fact that $V(t) = y^2(t)$ implies that:

$$|y(t)| \leq \max\left\{ \sqrt{2}\exp(-t)|y_0|, \varepsilon\|x\| \right\}, \text{ for all } t \in [0, t_{\max}) \tag{3.21}$$

Clearly, inequality (3.21) implies that as long as the solution of (3.15) exists, $y(t)$ is bounded. A standard contradiction argument in conjunction with the "Boundedness-Implies-Continuation" property for (3.15) (see Proposition 2.5 in [21]), implies that $t_{\max} = +\infty$. Inequality (3.16) is a direct consequence of inequality (3.21). The proof is complete. ◁

We select $M > 0$ sufficiently large and $r > 0$ sufficiently small such that inequality (3.16) holds with $\varepsilon < 1/\sqrt{2}$ and $a = 1 + \sqrt{2}/2$. The solution of the closed-loop system (4.1a), (4.2), where $k$ is defined by (4.5), exists for all $t \geq 0$. The existence of the solution is guaranteed by the following claim.

**Claim 2:** *For every $M > 0$, $r > 0$ and $(y_0, x_0, w) \in \Re \times \Re^n \times L_{loc}^\infty(\Re^+; \Re^+)$, the solution of (3.13), with initial condition $(x(0), y(0)) = (x_0, y_0)$ corresponding to input $w \in L_{loc}^\infty(\Re^+; \Re^+)$ exists for all $t \geq 0$. Moreover, for $M > 0$ sufficiently large and $r > 0$ sufficiently small, there exist $g \in K_\infty$ and $\xi \in \pi$ such that*

$$|x(t)| + |y(t)| \leq g(|(x_0, y_0)|), \text{ for all } t \in [0, \xi] \qquad (3.22)$$

$$|x(t)| \leq a, \text{ for all } t \geq \xi \qquad (3.23)$$

$$\xi \leq 1 + r + g(|x_0|) \qquad (3.24)$$

*where $a = 1 + \sqrt{2}/2$.*

**Proof of Claim 2:** Let $M > 0$, $r > 0$ and $(y_0, x_0, w) \in \Re \times \Re^n \times L_{loc}^\infty(\Re^+; \Re^+)$ be arbitrary. Consider a solution $(x(t), y(t))$ of (3.13) corresponding to arbitrary $w \in L_{loc}^\infty(\Re^+; \Re^+)$ with initial condition $(x(0), y(0)) = (x_0, y_0)$. By virtue of Proposition 2.5 in [21], there exists a maximal existence time for the solution denoted by $t_{\max} \leq +\infty$. Moreover, let $\pi := \{\tau_0, \tau_1, ...\}$ the set of sampling times (which may be finite if $t_{\max} < +\infty$) and $mp(t) := \max\{\tau \in \pi : \tau \leq t\}$. By virtue of (3.14) we have for all $t \in [0, t_{\max})$

$$|x(t)| \leq \max\{\sigma(|x_0|, 0), a\} \qquad (3.25)$$

Define

$$P := \max\{|f(x)| : |x| \leq \max\{\sigma(|x_0|, 0), a\}\} \text{ and } Q := \max\{|g(x)| : |x| \leq \max\{\sigma(|x_0|, 0), a\}\} \qquad (3.26)$$

Using (3.15), (3.26) in conjunction with Gronwall-Bellman's lemma we obtain the following inequality for all $t \in [0, t_{\max})$:

$$|y(t)| \leq |y(\tau)| \exp((M + 2Q)(t - \tau)) + P(t - \tau) \exp(Q(t - \tau)) \qquad (3.27)$$

where $\tau = mp(t)$. Using (3.27) and by induction we can show the following inequality for all $\tau_i \in \pi$:

$$|y(\tau_i)| \leq |y_0| \exp((M + 2Q)\tau_i) + P\tau_i \exp(Qr) \exp((M + 2Q)\tau_i) \qquad (3.28)$$

where we have used the fact that $\tau_{i+1} - \tau_i \leq r$. Estimate (3.27) in conjunction with (3.28) gives for all $t \in [0, t_{\max})$:

$$|y(t)| \leq [|y_0| + Pt \exp(Qr)] \exp((M + 2Q)t) \qquad (3.29)$$

A standard contradiction argument in conjunction with the "Boundedness-Implies-Continuation" property for (3.13) (see Proposition 2.5 in [21]), implies that $t_{\max} = +\infty$.

The existence of $\xi \in \pi$ such that (3.23) holds is a direct consequence of (3.14) and definitions $\gamma(s) := \dfrac{\sqrt{2}\,s}{\sqrt{1 + 4s^2}}$, $a = 1 + \sqrt{2}/2$. By virtue of (4.3) and Proposition 7 in [39] there exists $\beta \in K_\infty$ such

$$\xi \leq 1 + r + \beta(|x_0|) \qquad (3.30)$$

Finally, let $M > 0$ sufficiently large and $r > 0$ sufficiently small so that (3.16) holds for $a = 1 + \sqrt{2}/2$ and $\varepsilon < 1/\sqrt{2}$. For $x_0 \in \Re^n$ with $\sigma(|x_0|, 0) \leq a$ we obtain from (3.16) and (3.14) for all $t \geq 0$:

$$|x(t)| + |y(t)| \leq (1 + \varepsilon + \sqrt{2})(\sigma(|x_0|, 0) + \sigma(|y_0|, 0)) \qquad (3.31)$$

Using (3.25), (3.26), (3.29), (3.30) and (3.31) we guarantee the existence of $\tilde{\beta} \in K_\infty$ such that

$$|x(t)| + |y(t)| \leq \tilde{\beta}(|(x_0, y_0)|), \text{ for all } t \in [0, \xi] \qquad (3.32)$$

The existence of $g \in K_\infty$ satisfying (3.22) and (3.24) is a direct consequence of (3.30) and (3.32). The proof is complete. ◁

The fact that the robust global stabilization problem for (3.13) with sampled-data feedback applied with zero order hold is solvable with $M > 0$ is sufficiently large and $r > 0$ is sufficiently small is a consequence from all the above and Theorem 2.5. Indeed, we apply Theorem 2.5 with $n = 2$, $V_1 = |x|$, $V_2 = |y|$, $L = |x| + |y|$, $H = (x, y)$, $S(t) := \{(x, y) \in \Re \times \Re : |x| \leq a\}$, $\gamma_{1,2}(s) := \sqrt{2} s$, $\gamma_{2,1}(s) := \varepsilon s$, $\gamma_{1,1} \equiv 0$, $\gamma_{2,2} \equiv 0$, $\zeta \equiv 0$, $g^u \equiv 0$, $\eta \equiv 0$, $\tilde{\eta} \equiv 0$, $p^u \equiv 0$, $c(t) = \tilde{c}(t) = \nu(t) = \mu(t) = \kappa(t) \equiv 1 + r$, $g \equiv 0$, $p(s, w) := s + w$, for appropriate $a, b \in K_\infty$, $\sigma \in KL$ and $q \in N_2$. All hypotheses (H1-4) are satisfied by using the above definitions and previous results. Therefore, we conclude that the closed-loop system (3.13) with $M > 0$ being sufficiently large and $r > 0$ being sufficiently small, satisfies the UISS property from the input $w$ with zero gain.

The reader should notice that alternative sampled-data feedback designs for system (3.12) applied with zero order hold and positive sampling rate can be obtained by using the results [34,35], which, however, achieve semiglobal and practical stabilization. It should be emphasized that the feedback design obtained by using the trajectory-based small-gain results of the present work guarantee global and asymptotic stabilization. Moreover, robustness to perturbations of the sampling schedule is guaranteed (that is the reason for introducing the input $w$ in the closed-loop system (3.13)). ◁

## 4. A Delayed Chemostat Model

In this section we study the robust global feedback stabilization problem for system (1.2) under (1.3). More specifically, in order to emphasize the fact that the mapping $p : C^0([-r, 0]; (0, S_i)) \to (0, +\infty)$ is unknown, we will consider the stabilization problem of the equilibrium point $(S_s, X_s) \in (0, S_i) \times (0, +\infty)$ satisfying (1.4) for the uncertain chemostat model

$$\dot{X}(t) = \left( \min_{t-r \leq \tau \leq t} \mu(S(\tau)) + d(t) \left( \max_{t-r \leq \tau \leq t} \mu(S(\tau)) - \min_{t-r \leq \tau \leq t} \mu(S(\tau)) \right) - D(t) - b \right) X(t)$$
$$\dot{S}(t) = D(t)(S_i - S(t)) - K(S(t))\mu(S(t))X(t) \quad (4.1)$$
$$X(t) \in (0, +\infty), S(t) \in (0, S_i), D(t) \geq 0, d(t) \in [0,1]$$

where $d(t) \in [0,1]$ is the uncertainty. We will assume that:

**(H)** *There exists $S^* < S_s$ such that $\mu(S) > b$ for all $S \in [S^*, S_i]$.*

Hypothesis (H) is automatically satisfied for the case of a monotone specific growth rate. Hypothesis (H) can be satisfied for non-monotone specific growth rates (e.g., Haldane or generalized Haldane growth expressions). By using the trajectory-based small-gain Theorem 2.6 we can prove the following theorem.

**Theorem 4.1:** *Let $a > 0$ be a constant that satisfies*

$$\min_{S^* \leq S \leq S_i} \mu(S) - b > aD_s \frac{S_s}{S_i} \quad (4.2)$$

*Then the locally Lipschitz delay-free feedback law:*

$$D(t) = \frac{K(S(t))\mu(S(t))X(t) + aD_s(S_s - \min(S(t), S_s))}{S_i - \min(S(t), S_s)} \quad (4.3)$$

*achieves robust global stabilization of the equilibrium point $(\{S(\theta) = S_s, \theta \in [-r, 0]\}, X_s) \in C^0([-r, 0]; (0, S_i)) \times (0, +\infty)$, for the uncertain chemostat model (4.1) under hypothesis (H).*

It should be noted that the change of coordinates:

$$X = X_s \exp(x) \quad , \quad S = \frac{S_i \exp(y)}{G + \exp(y)} \tag{4.4a}$$

where $G := \frac{S_i}{S_s} - 1$ and the input transformation

$$D = D_s \exp(u) \tag{4.4b}$$

maps the set $(0, S_i) \times (0, +\infty)$ onto $\Re^2$ and the equilibrium point $(\{S(\theta) = S_s, \theta \in [-r,0]\}, X_s) \in C^0([-r,0]; (0, S_i)) \times (0, +\infty)$ of system (4.1) to the equilibrium point $0 \in \Re \times C^0([-r,0]; \Re)$ of the transformed control system:

$$\begin{aligned}
\dot{x}(t) &= \min_{t-r \leq \tau \leq t} \tilde{\mu}(y(\tau)) + d(t)\left( \max_{t-r \leq \tau \leq t} \tilde{\mu}(y(\tau)) - \min_{t-r \leq \tau \leq t} \tilde{\mu}(y(\tau)) \right) - D_s \exp(u(t)) - b \\
\dot{y}(t) &= D_s (G \exp(-y(t)) + 1)[\exp(u(t)) - (G + \exp(y(t))) g(y(t)) \exp(x(t))] \\
(x, y) &\in \Re^2, u(t) \in \Re, d(t) \in [0,1]
\end{aligned} \tag{4.5}$$

where

$$\begin{aligned}
\tilde{\mu}(y) &:= \mu\left( \frac{S_i \exp(y)}{G + \exp(y)} \right) \\
g(y) &:= \frac{X_s}{D_s S_i G} K\left( \frac{S_i \exp(y)}{G + \exp(y)} \right) \mu\left( \frac{S_i \exp(y)}{G + \exp(y)} \right)
\end{aligned} \tag{4.6}$$

In the new coordinates the feedback law (4.3) takes the form:

$$u(t) = \ln\left( g(y(t)) \exp(x(t)) \min(G + \exp(y(t)), G+1) + \frac{a}{G+1} \max(1 - \exp(y(t)), 0) \right) \tag{4.7}$$

The feedback law (4.7) (or (4.3)) is a delay-free feedback, which achieves global stabilization of $0 \in \Re \times C^0([-r,0]; \Re)$ for system (4.5) no matter how large the delay is. Furthermore, no knowledge of the maximum delay $r \geq 0$ is needed for the implementation of (4.7).

The proof of Theorem 4.1 is therefore equivalent to the proof of robust global asymptotic stability of the equilibrium point $0 \in \Re \times C^0([-r,0]; \Re)$ for system (4.5).

Before we give the proof of Theorem 4.1, it is important to understand the intuition that leads to the construction of the feedback law (4.7) and the ideas behind the proof of Theorem 4.1. To explain the procedure we follow the following arguments:

1. For the stabilization of the equilibrium point $0 \in \Re \times C^0([-r,0]; \Re)$, we first start with the stabilization of subsystem $\dot{y}(t) = D_s(G \exp(-y(t)) + 1)[\exp(u(t)) - (G + \exp(y(t))) g(y(t)) \exp(x(t))]$ with $x$ as input. Any feedback law which satisfies $u(t) = \ln(g(y(t)) \exp(x(t))(G+1))$ for $y(t) \geq 0$ and $u(t) > \ln(g(y(t)) \exp(x(t))(G + \exp(y(t))))$ for $y(t) < 0$ achieves ISS stabilization of the subsystem with $x$ as input.

2. In order to prove URGAS for the composite system by means of small-gain arguments one has to show the ISS property of the $x$-subsystem $\dot{x}(t) = \min_{t-r \leq \tau \leq t} \tilde{\mu}(y(\tau)) + d(t)\left( \max_{t-r \leq \tau \leq t} \tilde{\mu}(y(\tau)) - \min_{t-r \leq \tau \leq t} \tilde{\mu}(y(\tau)) \right) - D_s \exp(u(t)) - b$ with $y$ as input. Notice that the feedback selection from previous step gives $\dot{x}(t) = \min_{t-r \leq \tau \leq t} \tilde{\mu}(y(\tau)) + d(t)\left( \max_{t-r \leq \tau \leq t} \tilde{\mu}(y(\tau)) - \min_{t-r \leq \tau \leq t} \tilde{\mu}(y(\tau)) \right) - D_s g(y(t)) \exp(x(t))(G+1) - b$ for $y(t) \geq 0$ and $\dot{x}(t) < \min_{t-r \leq \tau \leq t} \tilde{\mu}(y(\tau)) + d(t) \left( \max_{t-r \leq \tau \leq t} \tilde{\mu}(y(\tau)) - \min_{t-r \leq \tau \leq t} \tilde{\mu}(y(\tau)) \right) - D_s g(y(t)) \exp(x(t))(G + \exp(y(t))) - b$ for $y(t) < 0$. The estimation of the derivative $\dot{x}(t)$ shows that the ISS inequality for the $x-$subsystem does not hold unless we have $\min_{t-r \leq \tau \leq t} \tilde{\mu}(y(\tau)) > b$ for all $t$ sufficiently large. By virtue of hypothesis (H) there exists $y^* < 0$, such that the ISS inequality for the the $x-$subsystem holds if $\min_{t-r \leq \tau \leq t} y(\tau) \geq y^*$ holds for all $t$ sufficiently large.

3. The feedback law $u(t) > \ln(g(y(t)) \exp(x(t))(G + \exp(y(t))))$ for $y(t) < 0$ is selected such that the inequality $\min_{t-r \leq \tau \leq t} y(\tau) \geq y^*$ holds for all initial conditions after a transient period. Since the ISS inequalities will hold only after this transient period the trajectory-based small-gain result Theorem 2.6 must be used for the proof of URGAS of the closed-loop system.

Schematically, we have:

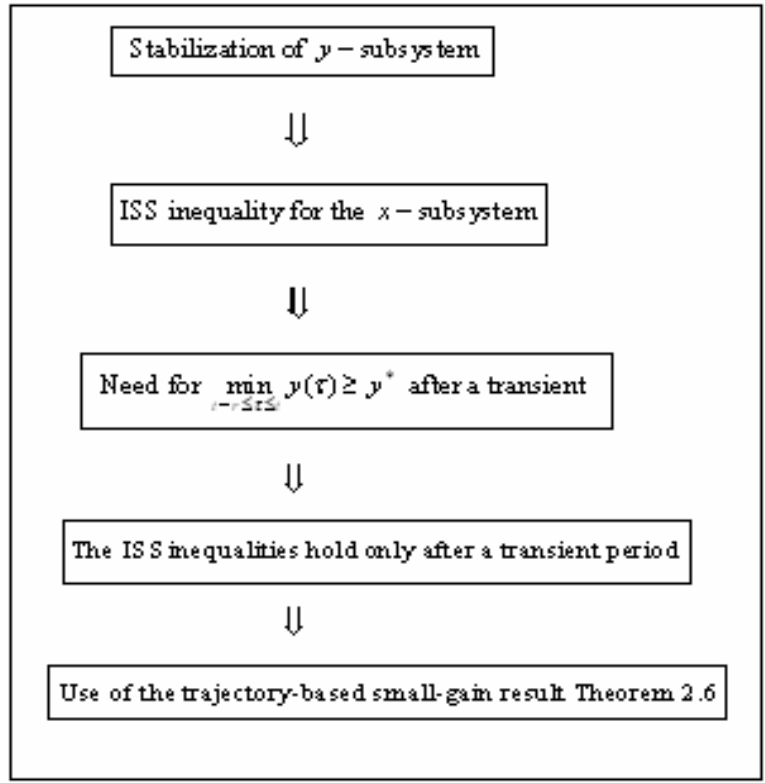

**Figure 1:** The intuition that leads to the construction of the feedback law (4.7) and the ideas behind the proof of Theorem 4.1

**Proof of Theorem 4.1:** Consider the solution $(x(t), y(t)) \in \Re^2$ of (4.5) with (4.7) with arbitrary initial condition $x(0) = x_0 \in \Re$, $T_r(0)y = y_0 \in C^0([-r,0]; \Re)$ and corresponding to arbitrary input $d \in M_D$. The following equations hold for system (4.5) with (4.7):

$$\dot{y}(t) = aD_s \frac{G \exp(-y(t)) + 1}{G+1} (1 - \exp(y(t))) \quad , \quad \text{if} \quad y(t) \leq 0$$
$$\dot{y}(t) = D_s g(y(t))(G \exp(-y(t)) + 1) \exp(x(t))(1 - \exp(y(t))) \quad , \quad \text{if} \quad y(t) > 0$$
(4.8)

Equations (4.8) imply that the function $V(t) = y^2(t)$ is non-increasing and consequently, we obtain:

$$|y(t)| \leq \|y_0\|_r, \text{ for all } t \in [0, t_{\max}) \tag{4.9}$$

Using the fact that $\mu : (0, S_i) \to (0, \mu_{\max}]$ and definition (4.6) of $\tilde{\mu}$, we get that $\tilde{\mu}(y) \leq \mu_{\max}$ for all $y \in \Re$. This implies the following differential inequality:

$$\dot{x}(t) \leq 2\mu_{\max} - b$$

which by direct integration yields the estimate:

$$x(t) \leq x_0 + (2\mu_{\max} - b)t, \text{ for all } t \in [0, t_{\max}) \tag{4.10}$$

Define $\kappa(s) := (G+1)D_s \max_{|y| \leq s} g(y)$. Inequalities (4.9), (4.10) imply that the following differential inequality holds:

$$\dot{x}(t) \geq -b - \frac{aD_s}{G+1} - \kappa(\|y_0\|_r)\exp(x_0 + (2\mu_{\max} - b)t)$$

which by direct integration yields the estimate:

$$x(t) \geq x_0 - \left(b + \frac{aD_s}{G+1}\right)t - \kappa(\|y_0\|_r)\exp(x_0)\frac{\exp((2\mu_{\max} - b)t) - 1}{2\mu_{\max} - b}, \text{ for all } t \in [0, t_{\max}) \tag{4.11}$$

Inequalities (4.9), (4.10), (4.11) and a standard contradiction argument show that system (4.5) with (4.7) is forward complete, i.e., $t_{\max} = +\infty$. Therefore, inequalities (4.9), (4.10), (4.11) hold for all $t \geq 0$ and since system (4.5) with (4.7) is autonomous, it follows that system (4.5) with (4.7) is Robustly Forward Complete (RFC, see Appendix B).

By considering (4.8) and the function $W(t) = \begin{cases} y^2(t) & \text{if } y(t) \leq 0 \\ 0 & \text{if } y(t) > 0 \end{cases}$, we obtain the existence of a positive definite function $\rho \in C^0(\Re_+; \Re_+)$ such that:

$$\dot{W}(t) \leq -\rho(W(t)), \text{ for all } t \geq 0 \tag{4.12}$$

Lemma 3.5 in [25] implies the existence of $\sigma \in KL$ such that for every $x(0) = x_0 \in \Re$, $T_r(0)y = y_0 \in C^0([-r,0];\Re)$ and $d \in M_D$ it holds that:

$$W(t) \leq \sigma(W(0), t), \text{ for all } t \geq 0 \tag{4.13}$$

Inequality (4.13) in conjunction with Proposition 7 in [39], shows the existence of $a \in K_\infty$ such that for every $x(0) = x_0 \in \Re$, $T_r(0)y = y_0 \in C^0([-r,0];\Re)$ and $d \in M_D$ there exists $\xi \geq r$ with $\xi \leq r + a(\|y_0\|_r)$ satisfying:

$$y(t-r) \geq -\frac{c}{2}, \text{ for all } t \geq \xi \tag{4.14}$$

where $c := \ln\left(\frac{S_i - S^*}{S^* G}\right) > 0$ and $S^* < S_s$ is the constant involved in hypothesis (H). Define $S := \Re \times C^0([-r,0];[-c,+\infty))$. Inequality (4.14) shows that $(x(t), T_r(t)y) \in S$ for all $t \geq \xi$ and that inequality (2.16) holds for appropriate $a \in K_\infty$ and $c(t) \equiv 1$.

Notice that for $(x(t), T_r(t)y) \in S$, the functionals

$$V_1(t) = \max_{\theta \in [-r,0]} \exp(2\sigma\theta)|z(t+\theta)|^2, \quad V_2 = |x(t)|^2 \tag{4.15}$$

where $\sigma > 0$ and

$$y(t) = c(\exp(z(t)) - 1) \tag{4.16}$$

are well-defined. Moreover, by considering the differential equations:

$$\dot{z}(t) = aD_s \frac{G\exp(c(1-\exp(z(t))))+1}{c(G+1)} \exp(-z(t))(1-\exp(c(\exp(z(t))-1))) \quad , \quad \text{if} \quad z(t) \leq 0$$

$$\dot{z}(t) = c^{-1}D_s g(c(\exp(z(t))-1))(G\exp(c(1-\exp(z(t))))+1)\exp(x(t)-z(t))(1-\exp(c(\exp(z(t))-1))) \quad , \quad \text{if} \quad z(t) > 0$$

we conclude from Lemma 3.5 in [25] that for every $\gamma_{1,2} \in K_\infty$ there exists $\sigma_1 \in KL$ such that:

$$V_1(t) \leq \max\left\{ \sigma_1(V_1(t_0), t-t_0), \sup_{t_0 \leq \tau \leq t} \gamma_{1,2}(V_2(\tau)) \right\}, \text{ for all } t \geq t_0 \geq 0 \tag{4.17}$$

Finally, using hypothesis (H) and definitions (4.15), we guarantee that there exists a positive definite function $\rho \in C^0(\Re_+; \Re_+)$ such that the following implication holds for every $\varepsilon > 0$:

"If $(1+\varepsilon)\ln\left( \dfrac{(G+\exp(c(\exp(\sqrt{V_1(t)})-1)))D_s \max_{|z|\leq \sqrt{V_1(t)}} g(c(\exp(z)-1))}{\min_{|z|\leq \exp(\sigma r)\sqrt{V_1(t)}} \tilde{\mu}(c(\exp(z)-1))-b - \dfrac{a}{G+1}D_s(1-\exp(c(\exp(-\sqrt{V_1(t)})-1)))} \right) \leq |x(t)|$ and

$(1+\varepsilon)\ln\left( \dfrac{\max_{|z|\leq \exp(\sigma r)\sqrt{V_1(t)}} \tilde{\mu}(c(\exp(z)-1))-b}{(G+\exp(c(\exp(-\sqrt{V_1(t)})-1)))D_s \min_{|z|\leq \sqrt{V_1(t)}} g(c(\exp(z)-1))} \right) \leq |x(t)|$ then

$2x(t)\dot{x}(t) \leq -\rho(x^2(t))$"

Therefore, Lemma 3.5 in [25] implies that there exists $\sigma_2 \in KL$ such that:

$$V_2(t) \leq \max\left\{ \sigma_2(V_2(t_0), t-t_0), \sup_{t_0 \leq \tau \leq t} \gamma_{2,1}(V_1(\tau)) \right\}, \text{ for all } t \geq t_0 \geq 0 \tag{4.18}$$

where

$$\gamma_{2,1}(s) := (1+\varepsilon)^2 (\ln(\max\{g_1(s), g_2(s)\}))^2$$

$$g_1(s) := \frac{(G+\exp(c(\exp(\sqrt{s})-1)))D_s \max_{|z|\leq \sqrt{s}} g(c(\exp(z)-1))}{\min_{|z|\leq \exp(\sigma r)\sqrt{s}} \tilde{\mu}(c(\exp(z)-1))-b - \dfrac{aD_s}{G+1}(1-\exp(c(\exp(-\sqrt{s})-1)))} \tag{4.19}$$

$$g_2(s) := \frac{\max_{|z|\leq \exp(\sigma r)\sqrt{s}} \tilde{\mu}(c(\exp(z)-1))-b}{(G+\exp(c(\exp(-\sqrt{s})-1)))D_s \min_{|z|\leq \sqrt{s}} g(c(\exp(z)-1))}$$

Inequalities (4.9), (4.10), (4.11), (4.17) and (4.18) guarantee that inequalities (2.10), (2.11), (2.13), (2.14), (2.15) and (2.17) hold for appropriate $\sigma \in KL$, $\nu \in K^+$, $a \in K_\infty$ with $c(t) \equiv 1$, $p \equiv 0$, $\gamma_{1,2}(s) := \gamma_{2,1}\left(\dfrac{s}{2}\right)$, $\gamma_{1,1}(s) = \gamma_{2,2}(s) \equiv 0$, $L := V_1 + V_2$ and $H(t, x, y) := \sqrt{x^2 + \|y\|_r^2}$. Finally, notice that the MAX-preserving mapping $\Gamma : \Re_+^2 \to \Re_+^2$ with $\Gamma_i(x) = \max_{j=1,2} \gamma_{i,j}(x_j)$ ($i = 1,2$) satisfies the cyclic small-gain conditions.

By virtue of Theorem 2.6 we conclude that the autonomous system (4.5) with (4.7) is URGAS. ◁

## 5. Conclusions

One of the most important obstacles in applying nonlinear small-gain results is the fact that the essential inequalities, which small-gain results utilize in order to prove stability properties, do not hold for all times: this feature excludes all available small-gain results from possible application. In this work novel small-gain results are presented, which can allow a transient period during which the solutions do not satisfy the usual inequalities required by previous small-gain results (Theorem 2.5 and Theorem 2.6). The obtained results allow the application of the small-gain methodology to various classes of systems which satisfy less demanding stability notions than the Input-to-Output Stability property.

The robust global feedback stabilization problem of an uncertain time-delay chemostat model is solved by means of the trajectory-based small-gain results. Future research will focus on the application of the trajectory-based small-gain results to Lotka-Volterra systems in Mathematical Biology (see [11]).

# Appendix A-Proofs of Theorem 2.5 and Theorem 2.6

**Proof of Theorem 2.5:** The proof is similar to the proof of Theorem 3.1 in [28] and consists of four steps:

<u>Step 1:</u> We show that for every $(t_0, x_0, u, d) \in \Re_+ \times \mathcal{X} \times M_U \times M_D$ the following inequality holds for all $t \in [\xi, t_{\max})$:

$$V(t) \leq MAX\left\{ Q(\mathbf{1}\sigma(L(\xi), 0)), Q\left(\mathbf{1}\zeta\left(\|u\|_{\mathcal{U}}\big]_{[\xi, t]}\right)\right) \right\} \qquad (A.1)$$

where $\xi \in \pi(t_0, x_0, u, d)$ is the time such that $\phi(t, t_0, x_0, u, d) \in S(t)$ for all $t \in [\xi, t_{\max})$ (recall Hypothesis (H2)).

This step is proved in exactly the same way as in the proof of Theorem 3.1 in [28] (using Fact V).

<u>Step 2:</u> We show that for every $(t_0, x_0, u, d) \in \Re_+ \times \mathcal{X} \times M_U \times M_D$, it holds that $t_{\max} = +\infty$.

The proof of Step 2 is exactly the same with the proof of Theorem 3.1 in [28]. The only difference is the additional use of inequality (2.6), which guarantees that the transition map is bounded during the transient period $t \in [t_0, \xi]$.

<u>Step 3:</u> We show that $\Sigma$ is RFC from the input $u \in M_U$.

Again the proof of Step 3 is exactly the same with the proof of Theorem 3.1 in [28]. The only difference is the additional use of inequalities (2.6) and (2.9).

<u>Step 4:</u> We prove the following claim.

**Claim:** *For every $\varepsilon > 0$, $k \in Z_+$, $R, T \geq 0$ there exists $\tau_k(\varepsilon, R, T) \geq 0$ such that for every $(t_0, x_0, u, d) \in \Re_+ \times \mathcal{X} \times M_U \times M_D$ with $t_0 \in [0, T]$ and $\|x_0\|_{\mathcal{X}} \leq R$ the following inequality holds:*

$$V(t) \leq MAX\left\{ Q(\mathbf{1}\varepsilon), \Gamma^{(k)}\big(Q(\mathbf{1}\sigma(L(\xi), 0))\big), G\left(\|u\|_{\mathcal{U}}\big]_{[t_0, t]}\right) \right\}, \text{ for all } t \geq \xi + \tau_k \qquad (A.2)$$

*Moreover, if $c \in K^+$ is bounded then for every $\varepsilon > 0$, $k \in Z_+$, $R \geq 0$ there exists $\tau_k(\varepsilon, R) \geq 0$ such that for every $(t_0, x_0, u, d) \in \Re_+ \times \mathcal{X} \times M_U \times M_D$ with $\|x_0\|_{\mathcal{X}} \leq R$ inequality (A.2) holds.*

<u>Proof of Step 4:</u> The proof of the claim is made by induction on $k \in Z_+$.

Inequality (A.2) for $k = 1$ is a direct consequence of inequalities (2.4), (2.9) and definition (2.12).

We notice that inequality (2.5) in conjunction with inequality (A.1) and Fact IV imply for all $t \geq \xi$:

$$L(t) \leq \max\left\{ \nu(t-t_0), c(t_0), a(\|x_0\|_X), p(Q(\mathbf{1}\sigma(L(\xi),0))), p\left(Q\left(\mathbf{1}\zeta\left(\|u\|_{\mathcal{U}}\right]_{[t_0,t]}\right)\right)\right), p^u\left(\|u\|_{\mathcal{U}}\right]_{[t_0,t]}\right) \right\} \quad (A.3)$$

Next suppose that for every $\varepsilon > 0$, $R, T \geq 0$ there exists $\tau_k(\varepsilon, R, T) \geq 0$ such that for every $(t_0, x_0, u, d) \in \Re_+ \times X \times M_U \times M_D$ with $t_0 \in [0,T]$ and $\|x_0\|_X \leq R$ (A.2) holds for some $k \in Z_+$. Let arbitrary $\varepsilon > 0$, $R, T \geq 0$, $(t_0, x_0, u, d) \in \Re_+ \times X \times M_U \times M_D$ with $t_0 \in [0,T]$ and $\|x_0\|_X \leq R$ be given. Notice that the weak semigroup property implies that $\pi(t_0, x_0, u, d) \cap [\xi + \tau_k, \xi + \tau_k + r] \neq \emptyset$. Let $t_k \in \pi(t_0, x_0, u, d) \cap [\xi + \tau_k, \xi + \tau_k + r]$. Then (2.4) implies:

$$V(t) \leq MAX\left\{ \mathbf{1}\sigma(L(t_k), t-t_k), \Gamma\left([V]_{[t_k,t]}\right), \mathbf{1}\zeta\left(\|u\|_{\mathcal{U}}\right]_{[t_k,t]}\right) \right\}, \text{ for all } t \geq t_k \quad (A.4)$$

Using inequalities (A.2), (A.3), (A.4), (2.9) and working in the same way as in the proof of Theorem 3.1 in [28] we can derive the following inequality:

$$V(t) \leq MAX\left\{ \mathbf{1}\sigma(L(t_k), t-\xi-\tau_k-r), Q(\mathbf{1}\varepsilon), \Gamma^{(k+1)}(Q(\mathbf{1}\sigma(L(\xi),0))), G\left(\|u\|_{\mathcal{U}}\right]_{[t_0,t]}\right) \right\},$$
$$\text{for all } t \geq \xi + \tau_k + r \quad (A.5)$$

Definition (2.12) in conjunction with (2.7), (2.9), inequality (A.3) and the facts that $t_k \leq \xi + \tau_k + r$, $t_0 \in [0,T]$ and $\|x_0\|_X \leq R$ implies that

$$\mathbf{1}\sigma(L(t_k), t-\xi-\tau_k-r) \leq MAX\left\{ \mathbf{1}\sigma(f(\varepsilon,T,R), t-\xi-\tau_k-r), G\left(\|u\|_{\mathcal{U}}\right]_{[t_0,t]}\right) \right\},$$
$$\text{for all } t \geq \xi + \tau_k + r \quad (A.6)$$

where

$$f(\varepsilon, T, R) := \max\left\{ \max_{0 \leq t \leq a(R) + C(T) + \tau_k(\varepsilon, R, T) + r} \nu(t), C(T), a(R), p(Q(\mathbf{1}\sigma(a(RC(T)),0))) \right\} \quad (A.7)$$

and

$$C(T) := \max_{0 \leq t \leq T} c(t) \quad (A.8)$$

The reader should notice that if $c \in K^+$ is bounded and $\tau_k$ is independent of $T$ then $f$ can be chosen to be independent of $T$ as well. The rest of proof of the claim follows from a combination of inequalities (A.5) and (A.6) and appropriate selection of $\tau_{k+1}$ (set $\tau_{k+1}(\varepsilon, R, T) = \tau_k(\varepsilon, R, T) + r + \tau(\varepsilon, R, T)$, where $\tau(\varepsilon, R, T) \geq 0$ satisfies $\sigma(f(\varepsilon, T, R), \tau) \leq \varepsilon$).

To finish the proof, let arbitrary $\varepsilon > 0$, $R, T \geq 0$, $(t_0, x_0, u, d) \in \Re_+ \times X \times M_U \times M_D$ and denote $Y(t) = H(t, \phi(t, t_0, x_0, u, d), u(t))$ for $t \geq t_0$. Using Fact IV, (2.11) and (A.1) we obtain for all $t \geq \xi$:

$$\|Y(t)\|_Y \leq \max\left\{ q(Q(\mathbf{1}\sigma(L(\xi),0))), q\left(Q\left(\mathbf{1}\zeta\left(\|u\|_{\mathcal{U}}\right]_{[\xi,t]}\right)\right)\right) \right\}$$

The above inequality in conjunction with (2.9) implies that

$$\|Y(t)\|_Y \leq \max\left\{ q(Q(\mathbf{1}\sigma(a(c(t_0)\|x_0\|_X),0))), q\left(Q\left(\mathbf{1}\sigma\left(g^u\left(\|u\|_{\mathcal{U}}\right]_{[t_0,t]}\right),0\right)\right)\right), q\left(Q\left(\mathbf{1}\zeta\left(\|u\|_{\mathcal{U}}\right]_{[t_0,t]}\right)\right)\right) \right\},$$
$$\text{for all } t \geq \xi \quad (A.9)$$

Using (2.8) and (A.9), we conclude that the following estimate holds for all $t \geq t_0$:

$$\|Y(t)\|_Y \leq \max \left\{ \begin{array}{l} q\big(Q\big(\mathbf{1}\sigma\big(a\big(c(t_0)\|x_0\|_X\big),0\big)\big)\big), a\big(c(t_0)\|x_0\|_X\big), \\ \eta\big(\|u\|_{\mathcal{U}[t_0,t]}\big), q\big(Q\big(\mathbf{1}\sigma\big(g^u\big(\|u\|_{\mathcal{U}[t_0,t]}\big),0\big)\big)\big), q\big(Q\big(\mathbf{1}\zeta\big(\|u\|_{\mathcal{U}[t_0,t]}\big)\big)\big) \end{array} \right\} \quad (A.10)$$

Inequality (A.10) shows that properties P1 and P2 of Lemma 2.16 in [23] hold for system $\Sigma$ with $V = \|H(t,x,u)\|_Y$ and $\gamma(s) := \max\{\eta(s), q(G(s))\}$. Moreover, if $c \in K^+$ is bounded then (A.10) implies that properties P1 and P2 of Lemma 2.17 in [23] hold for system $\Sigma$ with $V = \|H(t,x,u)\|_Y$ and $\gamma(s) := \max\{\eta(s), q(G(s))\}$.

Inequality (A.2) in conjunction with Fact III, (2.9), (A.8) and definition (2.12) guarantees that for every $\varepsilon > 0$, $k \in Z_+$, $R, T \geq 0$ there exists $\tau_k(\varepsilon, R, T) \geq 0$ such that for every $(t_0, x_0, u, d) \in \Re_+ \times X \times M_U \times M_D$ with $t_0 \in [0, T]$ and $\|x_0\|_X \leq R$ the following inequality holds:

$$V(t) \leq MAX\left\{ Q(\mathbf{1}\varepsilon), \Gamma^{(k)}\big(Q\big(\mathbf{1}\sigma\big(a(RC(T)),0\big)\big)\big), G\big(\|u\|_{\mathcal{U}[t_0,t]}\big) \right\}, \text{ for all } t \geq \xi + \tau_k \quad (A.11)$$

Notice that Fact I guarantees the existence of $k(\varepsilon, T, R) \in Z_+$ such that $Q(\mathbf{1}\varepsilon) \geq \Gamma^{(l)}\big(Q\big(\mathbf{1}\sigma\big(a(RC(T)),0\big)\big)\big)$ for all $l \geq k$. If $c \in K^+$ is bounded then $k$ is independent of $T$. Therefore by virtue of (A.11) and (2.7), property P3 of Lemma 2.16 in [23] holds for system $\Sigma$ with $V = \|H(t,x,u)\|_Y$ and $\gamma(s) := \max\{\eta(s), q(G(s))\}$. Moreover, if $c \in K^+$ is bounded then (A.11) and (2.7) imply that property P3 of Lemma 2.17 in [23] hold for system $\Sigma$ with $V = \|H(t,x,u)\|_Y$ and $\gamma(s) := \max\{\eta(s), q(G(s))\}$.

The proof of Theorem 2.4 is thus completed with the help of Lemma 2.16 (or Lemma 2.17) in [23]. ◁

**Proof of Theorem 2.6:** By virtue of Lemma 3.3 in [20] we have to show that $\Sigma$ is Robustly Forward Complete (RFC) and satisfies the *Robust Output Attractivity Property*, i.e. for every $\varepsilon > 0$, $T \geq 0$ and $R \geq 0$, there exists a $\tau := \tau(\varepsilon, T, R) \geq 0$, such that:

$$\|x_0\|_X \leq R, t_0 \in [0, T] \Rightarrow \|H(t, \phi(t, t_0, x_0, u_0, d), 0)\|_Y \leq \varepsilon, \forall t \in [t_0 + \tau, +\infty), \forall d \in M_D$$

The reader should notice that Lemma 3.3 in [20] assumes the classical semigroup property; however the semigroup property is not used in the proof of Lemma 3.3 in [20]. Consequently, Lemma 3.3 in [20] holds as well for systems satisfying the weak semigroup property.

Moreover, Lemma 3.5 in [21] guarantees that system $\Sigma$ is URGAOS in case that $\Sigma := (X, Y, M_U, M_D, \phi, \pi, H)$ is $T-$ periodic for certain $T > 0$.

Again the proof consists of four steps:

<u>Step 1:</u> We show that for every $(t_0, x_0, d) \in \Re_+ \times X \times M_D$ the following inequality holds for all $t \in [\xi, t_{\max})$:

$$V(t) \leq Q(\mathbf{1}\sigma(L(\xi), 0)) \quad (A.12)$$

where $\xi \in \pi(t_0, x_0, u_0, d)$ is the time such that $\phi(t, t_0, x_0, u_0, d) \in S(t)$ for all $t \in [\xi, t_{\max})$ (recall Hypothesis (H2)).

<u>Step 2:</u> We show that for every $(t_0, x_0, d) \in \Re_+ \times X \times M_D$, it holds that $t_{\max} = +\infty$.

<u>Step 3:</u> We show that $\Sigma$ is RFC.

<u>Step 4:</u> We prove the following claim.

**Claim:** *For every* $\varepsilon > 0$, $k \in Z_+$, $R, T \geq 0$ *there exists* $\tau_k(\varepsilon, R, T) \geq 0$ *such that for every* $(t_0, x_0, d) \in \Re_+ \times \mathcal{X} \times M_D$ *with* $t_0 \in [0, T]$ *and* $\|x_0\|_\mathcal{X} \leq R$ *the following inequality holds:*

$$V(t) \leq MAX\left\{Q(\mathbf{1}\varepsilon), \Gamma^{(k)}\left(Q(\mathbf{1}\sigma(L(\xi),0))\right)\right\}, \text{ for all } t \geq \xi + \tau_k \tag{A.13}$$

The proofs of the above steps are exactly the same with the proof of Theorem 2.5 and are omitted. The difference between inequalities (2.9) and (2.17) does not play any role. To finish the proof, let arbitrary $\varepsilon > 0$, $R, T \geq 0$, $(t_0, x_0, d) \in \Re_+ \times \mathcal{X} \times M_D$ and denote $Y(t) = H(t, \phi(t, t_0, x_0, u_0, d), 0)$ for $t \geq t_0$.

Inequality (A.13) in conjunction with (2.17) and (A.8) guarantees that for every $\varepsilon > 0$, $k \in Z_+$, $R, T \geq 0$ there exists $\tau_k(\varepsilon, R, T) \geq 0$ such that for every $(t_0, x_0, d) \in \Re_+ \times \mathcal{X} \times M_D$ with $t_0 \in [0, T]$ and $\|x_0\|_\mathcal{X} \leq R$ the following inequality holds:

$$V(t) \leq MAX\left\{Q(\mathbf{1}\varepsilon), \Gamma^{(k)}\left(Q(\mathbf{1}\sigma(a(R) + C(T), 0))\right)\right\}, \text{ for all } t \geq \xi + \tau_k \tag{A.14}$$

Notice that Fact I guarantees the existence of $k(\varepsilon, T, R) \in Z_+$ such that $Q(\mathbf{1}\varepsilon) \geq \Gamma^{(l)}\left(Q(\mathbf{1}\sigma(a(R) + C(T), 0))\right)$ for all $l \geq k$. Therefore, (A.14) implies that for every $\varepsilon > 0$, $R, T \geq 0$, there exists $\tau(\varepsilon, R, T) \geq 0$ such that for every $(t_0, x_0, d) \in \Re_+ \times \mathcal{X} \times M_D$ with $t_0 \in [0, T]$ and $\|x_0\|_\mathcal{X} \leq R$, it holds that

$$V(t) \leq Q(\mathbf{1}\varepsilon), \text{ for all } t \geq \xi + \tau \tag{A.15}$$

It follows from inequalities (2.11) and (A.15) that for every $\varepsilon > 0$, $R, T \geq 0$, there exists $\tau(\varepsilon, R, T) \geq 0$ such that for every $(t_0, x_0, d) \in \Re_+ \times \mathcal{X} \times M_D$ with $t_0 \in [0, T]$ and $\|x_0\|_\mathcal{X} \leq R$, it holds that

$$\|Y(t)\|_\mathcal{Y} \leq q(Q(\mathbf{1}\varepsilon)), \text{ for all } t \geq \xi + \tau \tag{A.16}$$

Therefore by virtue of (A.16) and (2.16), the *Robust Output Attractivity Property* holds for system $\Sigma$. The proof is complete. ◁

## Appendix B-Basic Notions

To make our work self-contained, we introduce some notions which are essential to the system theoretic framework presented in [21,22,23]. The abstract system theoretic framework used in [21,22,23] is utilized in the present work.

**The notion of a Control System-Definition 2.1 in [23]:** *A control system* $\Sigma := (\mathcal{X}, \mathcal{Y}, M_U, M_D, \phi, \pi, H)$ *with outputs consists of*

**(i)** *a set* $U$ *(control set) which is a subset of a normed linear space* $\mathcal{U}$ *with* $0 \in U$ *and a set* $M_U \subseteq \mathcal{M}(U)$ *(allowable control inputs) which contains at least the identically zero input* $u_0$,

**(ii)** *a set* $D$ *(disturbance set) and a set* $M_D \subseteq \mathcal{M}(D)$, *which is called the "set of allowable disturbances",*

**(iii)** *a pair of normed linear spaces* $\mathcal{X}, \mathcal{Y}$ *called the "state space" and the "output space", respectively,*

**(iv)** *a continuous map* $H : \Re_+ \times \mathcal{X} \times U \to \mathcal{Y}$ *that maps bounded sets of* $\Re_+ \times \mathcal{X} \times U$ *into bounded sets of* $\mathcal{Y}$, *called the "output map",*

**(v)** *a set-valued map* $\Re_+ \times \mathcal{X} \times M_U \times M_D \ni (t_0, x_0, u, d) \to \pi(t_0, x_0, u, d) \subseteq [t_0, +\infty)$, *with* $t_0 \in \pi(t_0, x_0, u, d)$ *for all* $(t_0, x_0, u, d) \in \Re_+ \times \mathcal{X} \times M_U \times M_D$, *called the set of "sampling times"*

**(vi)** *and the map* $\phi : A_\phi \to \mathcal{X}$ *where* $A_\phi \subseteq \Re_+ \times \Re_+ \times \mathcal{X} \times M_U \times M_D$, *called the "transition map", which has the following properties:*

**1) Existence:** *For each $(t_0, x_0, u, d) \in \Re_+ \times X \times M_U \times M_D$, there exists $t > t_0$ such that $[t_0, t] \times (t_0, x_0, u, d) \subseteq A_\phi$.*

**2) Identity Property:** *For each $(t_0, x_0, u, d) \in \Re_+ \times X \times M_U \times M_D$, it holds that $\phi(t_0, t_0, x_0, u, d) = x_0$.*

**3) Causality:** *For each $(t, t_0, x_0, u, d) \in A_\phi$ with $t > t_0$ and for each $(\tilde{u}, \tilde{d}) \in M_U \times M_D$ with $(\tilde{u}(\tau), \tilde{d}(\tau)) = (u(\tau), d(\tau))$ for all $\tau \in [t_0, t]$, it holds that $(t, t_0, x_0, \tilde{u}, \tilde{d}) \in A_\phi$ with $\phi(t, t_0, x_0, u, d) = \phi(t, t_0, x_0, \tilde{u}, \tilde{d})$.*

**4) Weak Semigroup Property:** *There exists a constant $r > 0$, such that for each $t \geq t_0$ with $(t, t_0, x_0, u, d) \in A_\phi$:*
(a) $(\tau, t_0, x_0, u, d) \in A_\phi$ for all $\tau \in [t_0, t]$,
(b) $\phi(t, \tau, \phi(\tau, t_0, x_0, u, d), u, d) = \phi(t, t_0, x_0, u, d)$ for all $\tau \in [t_0, t] \cap \pi(t_0, x_0, u, d)$,
(c) if $(t+r, t_0, x_0, u, d) \in A_\phi$, then it holds that $\pi(t_0, x_0, u, d) \cap [t, t+r] \neq \emptyset$.
(d) for all $\tau \in \pi(t_0, x_0, u, d)$ with $(\tau, t_0, x_0, u, d) \in A_\phi$ we have $\pi(\tau, \phi(\tau, t_0, x_0, u, d), u, d) = \pi(t_0, x_0, u, d) \cap [\tau, +\infty)$.

**The BIC and RFC properties-Definition 2.4 in [23]:** *Consider a control system $\Sigma := (X, Y, M_U, M_D, \phi, \pi, H)$ with outputs. We say that system $\Sigma$*

(i) *has the **"Boundedness-Implies-Continuation" (BIC)** property if for each $(t_0, x_0, u, d) \in \Re_+ \times X \times M_U \times M_D$, there exists a maximal existence time, i.e., there exists $t_{\max} := t_{\max}(t_0, x_0, u, d) \in (t_0, +\infty]$, such that $A_\phi = \bigcup_{(t_0, x_0, u, d) \in \Re_+ \times X \times M_U \times M_D} [t_0, t_{\max}) \times \{(t_0, x_0, u, d)\}$. In addition, if $t_{\max} < +\infty$ then for every $M > 0$ there exists $t \in [t_0, t_{\max})$ with $\|\phi(t, t_0, x_0, u, d)\|_X > M$.*

(ii) *is **robustly forward complete (RFC) from the input** $u \in M_U$ if it has the BIC property and for every $r \geq 0$, $T \geq 0$, it holds that*

$$\sup\left\{\|\phi(t_0 + s, t_0, x_0, u, d)\|_X \; ; u \in M(B_U[0, r]) \cap M_U \; , s \in [0, T], \|x_0\|_X \leq r, t_0 \in [0, T], d \in M_D \right\} < +\infty$$

**The notion of a robust equilibrium point-Definition 2.5 in [23]:** *Consider a control system $\Sigma := (X, Y, M_U, M_D, \phi, \pi, H)$ and suppose that $H(t, 0, 0) = 0$ for all $t \geq 0$. We say that $0 \in X$ is a **robust equilibrium point from the input** $u \in M_U$ for $\Sigma$ if*

(i) *for every $(t, t_0, d) \in \Re_+ \times \Re_+ \times M_D$ with $t \geq t_0$ it holds that $\phi(t, t_0, 0, u_0, d) = 0$.*

(ii) *for every $\varepsilon > 0$, $T, h \in \Re_+$ there exists $\delta := \delta(\varepsilon, T, h) > 0$ such that for all $(t_0, x, u) \in [0, T] \times X \times M_U$, $\tau \in [t_0, t_0 + h]$ with $\|x\|_X + \sup_{t \geq 0}\|u(t)\|_U < \delta$ it holds that $(\tau, t_0, x, u, d) \in A_\phi$ for all $d \in M_D$ and*

$$\sup\left\{\|\phi(\tau, t_0, x, u, d)\|_X \; ; d \in M_D \; , \tau \in [t_0, t_0 + h], t_0 \in [0, T] \right\} < \varepsilon$$

Next we present the Input-to-Output Stability property for the class of systems described previously (see also [14,40] for finite-dimensional, time-invariant dynamic systems).

**The notions of IOS, UIOS, ISS and UISS-Definition 2.14 in [23]:** *Consider a control system $\Sigma := (X, Y, M_U, M_D, \phi, \pi, H)$ with outputs and the BIC property and for which $0 \in X$ is a robust equilibrium point from the input $u \in M_U$. Suppose that $\Sigma$ is RFC from the input $u \in M_U$. If there exist functions $\sigma \in KL$, $\beta \in K^+$, $\gamma \in N_1$ such that the following estimate holds for all $u \in M_U$, $(t_0, x_0, d) \in \Re_+ \times X \times M_D$ and $t \geq t_0$:*

$$\|H(t, \phi(t, t_0, x_0, u, d), u(t))\|_Y \leq \sigma\big(\beta(t_0)\|x_0\|_X, t - t_0\big) + \sup_{t_0 \leq \tau \leq t} \gamma\big(\|u(\tau)\|_U\big)$$

*then we say that $\Sigma$ satisfies the **Input-to-Output Stability (IOS)** property from the input $u \in M_U$ with gain $\gamma \in N$. Moreover, if $\beta \in K^+$ may be chosen as $\beta(t) \equiv 1$ then we say that $\Sigma$ satisfies the Uniform Input-to-Output Stability (UIOS) property from the input $u \in M_U$ with gain $\gamma \in N_1$.*

*For the special case of the identity output mapping, i.e., $H(t,x,u) := x$, the (Uniform) Input-to-Output Stability property from the input $u \in M_U$ is called (Uniform) Input-to-State Stability ((U) ISS) property from the input $u \in M_U$. When $U = \{0\}$ (the no-input case) and $\Sigma$ satisfies the (U)IOS property, then we say that $\Sigma$ satisfies the (Uniform) Robust Global Asymptotic Output Stability (RGAOS) property. When $U = \{0\}$ (the no-input case) and $\Sigma$ satisfies the (Uniform) ISS property, then we say that $\Sigma$ satisfies the (Uniform) Robust Global Asymptotic Stability (RGAS) property.*

## References

[1] Angeli, D., P. De Leenheer and E. D. Sontag, "A Small-Gain Theorem for Almost Global Convergence of Monotone Systems", *Systems and Control Letters*, 52(5), 2004, 407-414.

[2] Angeli, D. and A. Astolfi, "A Tight Small-Gain Theorem for not necessarily ISS Systems", *Systems and Control Letters*, 56, 2007, 87-91.

[3] Antonelli, R. and A. Astolfi, "Nonlinear Controller Design for Robust Stabilization of Continuous Biological Reactors", *Proceedings of the IEEE Conference on Control Applications*, Anchorage, AL, September 2000.

[4] Chen, Z. and J. Huang, "A simplified small gain theorem for time-varying nonlinear systems", *IEEE Transactions on Automatic Control*, 50(11), 2005, 1904-1908.

[5] Dashkovskiy, S., B. S. Ruffer and F. R. Wirth, "An ISS Small-Gain Theorem for General Networks", *Mathematics of Control, Signals and Systems*, 19, 2007, 93-122.

[6] De Leenheer, P. and H. L. Smith, "Feedback Control for Chemostat Models", *Journal of Mathematical Biology*, 46, 2003, 48-70.

[7] Enciso, G. A. and E. D. Sontag, "Global Attractivity, I/O Monotone Small-Gain Theorems, and Biological Delay Systems", *Discrete and Continuous Dynamical Systems*, 14(3), 2006, 549-578.

[8] Gouze, J. L. and G. Robledo, "Robust Control for an Uncertain Chemostat Model", *International Journal of Robust and Nonlinear Control*, 16(3), 2006, 133-155.

[9] Grune, L., "Input-to-State Dynamical Stability and its Lyapunov Function Characterization", *IEEE Transactions on Automatic Control*, 47(9), 2002, 1499-1504.

[10] Harmard, J., A. Rapaport and F. Mazenc, "Output tracking of continuous bioreactors through recirculation and by-pass", *Automatica*, 42, 2006, 1025-1032.

[11] Hofbauer, J. and K. Sigmund, *Evolutionary Games and Population Dynamics*, Cambridge University Press, 2002.

[12] Ito, H., "Stability Criteria for Interconnected iISS and ISS Systems Using Scaling of Supply Rates", *Proceedings of the American Control Conference*, 2, 2004, 1055-1060.

[13] Ito, H. and Z.-P. Jiang, "Necessary and Sufficient Small-Gain Conditions for Integral Input-to-State Stable Systems: A Lyapunov Perspective", *IEEE Transactions on Automatic Control*, 54(10), 2009, 2389-2404.

[14] Jiang, Z.P., A. Teel and L. Praly, "Small-Gain Theorem for ISS Systems and Applications", *Mathematics of Control, Signals and Systems*, 7, 1994, 95-120.

[15] Jiang, Z.P., I.M.Y. Mareels and Y. Wang, "A Lyapunov Formulation of the Nonlinear Small-Gain Theorem for Interconnected Systems", *Automatica*, 32, 1996, 1211-1214.

[16] Jiang, Z.P. and I.M.Y. Mareels, "A Small-Gain Control Method for Nonlinear Cascaded Systems with Dynamic Uncertainties", *IEEE Transactions on Automatic Control*, 42, 1997, 292-308.

[17] Jiang, Z.P., Y. Lin and Y. Wang, "Nonlinear Small-Gain Theorems for Discrete-Time Feedback Systems and Applications", *Automatica*, 40(12), 2004, 2129-2134.

[18] Jiang, Z.P. and Y. Wang, "A Generalization of the Nonlinear Small-Gain Theorem for Large-Scale Complex Systems", *Proceedings of the 7th World Congress on Intelligent Control and Automation*, Chongqing, China, 2008, 1188-1193.

[19] Karafyllis, I. and J. Tsinias, "Non-Uniform in Time ISS and the Small-Gain Theorem", *IEEE Transactions on Automatic Control*, 49(2), 2004, 196-214.

[20] Karafyllis, I., "The Non-Uniform in Time Small-Gain Theorem for a Wide Class of Control Systems with Outputs", *European Journal of Control*, 10(4), 2004, 307-322.

[21] Karafyllis, I., "A System-Theoretic Framework for a Wide Class of Systems I: Applications to Numerical Analysis", *Journal of Mathematical Analysis and Applications*, 328(2), 2007, 876-899.

[22] Karafyllis, I., "A System-Theoretic Framework for a Wide Class of Systems II: Input-to-Output Stability", *Journal of Mathematical Analysis and Applications*, 328(1), 2007, 466-484.

[23] Karafyllis, I. and Z.-P. Jiang, "A Small-Gain Theorem for a Wide Class of Feedback Systems with Control Applications", *SIAM Journal Control and Optimization*, 46(4), 2007, 1483-1517.

[24] Karafyllis, I., C. Kravaris, L. Syrou and G. Lyberatos, "A Vector Lyapunov Function Characterization of Input-to-State Stability with Application to Robust Global Stabilization of the Chemostat", *European Journal of Control*, 14(1), 2008, 47-61.

[25] Karafyllis, I. and C. Kravaris, "Global Stability Results for Systems under Sampled-Data Control", *International Journal of Robust and Nonlinear Control*, 19, 2009, 1105-1128.
[26] Karafyllis, I., and C. Kravaris, "Robust Global Stabilizability by Means of Sampled-Data Control with Positive Sampling Rate", *International Journal of Control*, 82(4), 2009, 755-772.
[27] Karafyllis, I., C. Kravaris and N. Kalogerakis, "Relaxed Lyapunov Criteria for Robust Global Stabilization of Nonlinear Systems", *International Journal of Control*, 82(11), 2009, 2077-2094.
[28] Karafyllis, I. and Z.-P. Jiang, "A Vector Small-Gain Theorem for General Nonlinear Control Systems", submitted to *IEEE Transactions on Automatic Control*. A short version was published in the *Proceedings of the 48th IEEE Conference on Decision and Control 2009*, Shanghai, China, 2009, pp. 7996-8001. A preliminary version is available at http://arxiv.org/abs/0904.0755.
[29] Mailleret, L. and O. Bernard, "A Simple Robust Controller to Stabilize an Anaerobic Digestion Process" in *Proceedings of the 8th Conference on Computer Applications in Biotechnology*, 2001, 213-218.
[30] Malisoff, M. and F. Mazenc, *Constructions of Strict Lyapunov Functions*, Springer Verlag, London, 2009.
[31] Mazenc, F., M. Malisoff, J. Harmand, "Stabilization and Robustness Analysis for a Chemostat Model with Two Species and Monod Growth Rates via a Lyapunov Approach", *Proceedings of the 46th IEEE Conference on Decision and Control*, New Orleans, 2007.
[32] Mazenc, F., M. Malisoff and P. De Leenheer, "On the Stability of Periodic Solutions in the Perturbed Chemostat", *Mathematical Biosciences and Engineering*, 4(2), 2007, 319-338.
[33] Mazenc, F., M. Malisoff and J. Harmand, "Further Results on Stabilization of Periodic Trajectories for a Chemostat with Two Species", *IEEE Transactions on Automatic Control*, 53(1), 2008, 66-74.
[34] Nesic, D. and A.R. Teel, "Sampled-Data Control of Nonlinear Systems: An Overview of Recent Results", Perspectives on Robust Control, R.S.O. Moheimani (Ed.), Springer-Verlag: New York, 2001, 221-239.
[35] Nesic, D. and A. R. Teel, "Stabilization of Sampled-Data Nonlinear Systems via Backstepping on their Euler Approximate Model", *Automatica*, 42, 2006, 1801-1808.
[36] Ruffer, B. S., *Monotone Dynamical Systems, Graphs and Stability of Large-Scale Interconnected Systems*, PhD Thesis, University of Bremen, Germany, 2007.
[37] Smith, H. and P. Waltman, *The Theory of the Chemostat. Dynamics of Microbial Competition*, Cambridge Studies in Mathematical Biology, 13, Cambridge University Press: Cambridge, 1995.
[38] Sontag, E.D., "Smooth Stabilization Implies Coprime Factorization", *IEEE Transactions on Automatic Control*, 34, 1989, 435-442.
[39] Sontag, E.D., "Comments on Integral Variants of ISS", *Systems and Control Letters*, 34, 1998, 93-100.
[40] Sontag, E.D. and Y. Wang, "Notions of Input to Output Stability", *Systems and Control Letters*, 38, 1999, 235-248.
[41] Sontag, E.D. and B. Ingalls, "A Small-Gain Theorem with Applications to Input/Output Systems, Incremental Stability, Detectability, and Interconnections*", Journal of the Franklin Institute*, 339, 2002, 211-229.
[42] Teel, A., "A Nonlinear Small Gain Theorem For the Analysis of Control Systems With Saturations", *IEEE Transactions on Automatic Control*, 41, 1996, 1256-1270.
[43] Teel, A.R., "Input-to-State Stability and the Nonlinear Small Gain Theorem", Preprint, 2005.
[44] Wang, L. and G. Wolkowicz, "A delayed chemostat model with general nonmonotone response functions and differential removal rates", *Journal of Mathematical Analysis and Applications*, 321, 2006, 452-468.
[45] Wolkowicz, G. and H. Xia, "Global Asymptotic Behavior of a chemostat model with discrete delays", *SIAM Journal on Applied Mathematics*, 57, 1997, 1019-1043.
[46] Zhu, L. and X. Huang, "Multiple limit cycles in a continuous culture vessel with variable yield", *Nonlinear Analysis*, 64, 2006, 887-894.
[47] Zhu, L., X. Huang and H. Su, "Bifurcation for a functional yield chemostat when one competitor produces a toxin", *Journal of Mathematical Analysis and Applications*, 329, 2007, 891-903.